# Almost-$C^1$ splines: Biquadratic splines on unstructured quadrilateral meshes and their application to fourth order problems

Thomas Takacs[a], Deepesh Toshniwal[b,*]

[a]Johannes Kepler University Linz, Austria
[b]Delft Institute of Applied Mathematics, Delft University of Technology, The Netherlands

**Abstract**

Isogeometric Analysis generalizes classical finite element analysis and intends to integrate it with the field of Computer-Aided Design. A central problem in achieving this objective is the reconstruction of analysis-suitable models from Computer-Aided Design models, which is in general a non-trivial and time-consuming task. In this article, we present a novel spline construction, that enables model reconstruction as well as simulation of high-order PDEs on the reconstructed models. The proposed *almost-$C^1$ splines* are biquadratic splines on fully unstructured quadrilateral meshes (without restrictions on placements or number of extraordinary vertices). They are $C^1$ smooth almost everywhere, that is, at all vertices and across most edges, and in addition almost (i.e. approximately) $C^1$ smooth across all other edges. Thus, the splines form $H^2$-nonconforming analysis-suitable discretization spaces. This is the lowest-degree unstructured spline construction that can be used to solve fourth-order problems. The associated spline basis is non-singular and has several B-spline-like properties (e.g., partition of unity, non-negativity, local support), the almost-$C^1$ splines are described in an explicit Bézier-extraction-based framework that can be easily implemented. Numerical tests suggest that the basis is well-conditioned and exhibits optimal approximation behavior.

## Contents

*Corresponding author
*Email addresses:* thomas.takacs@jku.at (Thomas Takacs), d.toshniwal@tudelft.nl (Deepesh Toshniwal)

# 1. Introduction

In this article we present a new approach for building analysis-suitable biquadratic spline spaces on fully unstructured quadrilateral meshes, so-called almost-$C^1$ splines. In particular, with the recent mixed smoothness spline construction from [64] as the starting point, we build almost $C^1$ smooth, $H^2$-nonconforming spline spaces that can be used to solve fourth-order problems. We test this on several model problems, such as the biharmonic problem, Kirchhoff–Love thin shells, a surface Cahn–Hilliard model and the surface Laplace–Beltrami eigenvalue problem. We obtain almost-$C^1$ splines by employing approximate $C^1$ smoothness at a fixed number of mesh edges (depending only on mesh topology and independent of its refinement level) and classical, parametric $C^1$ smoothness at all other edges. The motivation for doing so is manifold.

## 1.1. Motivation

This spline construction is inspired by the general aim of Isogeometric Analysis (*IGA*), introduced in [24], which is the integration of a finite element-like analysis within Computer-Aided Design (*CAD*) [17], cf. [23]. Achieving this objective would lead to a uniform and significantly more efficient design-through-analysis workflow for many engineering applications. In contrast, in the current setup a significant portion of the time is spent neither on design nor on analysis but is dominated by generating analysis-suitable reconstructions of CAD models (studies suggest up to 80%, cf. [6]).

Thus, our focus lies on developing splines that enable model reconstruction as well as simulation of high-order PDEs on the reconstructed models (cf. Figure 1). The goal is to achieve such analysis-suitable spline reconstructions for arbitrary topology geometries. While IGA was first introduced on single NURBS patches, that is, on structured quadrilateral meshes, to handle general bivariate geometries (planar domains and surfaces), one must be able to define analysis-suitable splines on general unstructured meshes.

To increase the geometric flexibility of the construction, the almost-$C^1$ splines are only approximately $C^1$ smooth near extraordinary vertices. Nonetheless, it is known that approximately $C^1$ smooth function spaces can be used to build optimally convergent finite element methods for partial differential equations (PDEs) in variational form requiring $H^2$ regularity; see, e.g., [54]. By relaxing the $C^1$ smoothness constraints in a specific manner, the singularities that appear for parametrically smooth spline constructions can be avoided (cf. [63]) and, potentially, better numerical convergence than non-singular geometrically smooth spline constructions (e.g. as in [42]) could be obtained.

Since we provide explicit descriptions of the almost-$C^1$ splines, they form a viable alternative to other approximately smooth constructions that are based on a weak imposition of smoothness, such as Nitsche's method [21, 40], the mortar method [4, 22, 39] or a mixed approach [46]. A similar, explicit construction for approximately $C^1$ smooth isogeometric multi-patch spaces is presented in [68].

Finally, the proposed almost-$C^1$ splines allow a construction with a highly local footprint and fewer restrictions than the available parametrically or geometrically smooth alternatives, cf. [25].

With this motivation in mind, we begin by highlighting some prior work in these areas in Section 1.2, then we summarize the main properties of almost-$C^1$ splines in Section 1.3 and present a short outline of this article in Section 1.4.

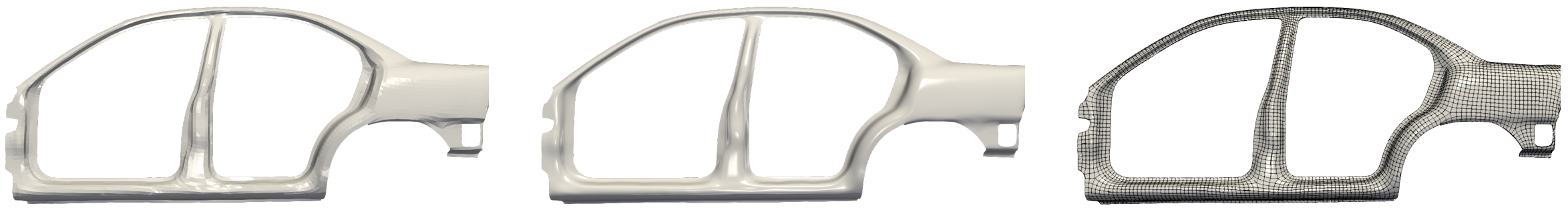

(a) watertight quad mesh (b) reconstructed spline surface (c) underlying Bézier mesh

Figure 1: The above figures show how almost-$C^1$ splines can enable analysis-suitable model reconstructions. For a crash simulation vehicle model of the Dodge Neon (not shown) [44], figure (a) shows a watertight bilinear quadrilateral mesh generated using Rhinoceros®. Figure (b) shows an analysis-suitable spline model reconstructed using almost-$C^1$ splines; figure (c) shows the underlying Bézier mesh. Note that watertight mesh generation for CAD models is itself an active area of research; for instance, see the recent article [53]. Note also that we used the quad-mesh in figure (a) to automatically compute control points used in figure (b); in practice this reconstruction step can be improved by using information from the original (CAD or finite element) model.

*1.2. Smooth splines over unstructured quadrilateral meshes*

In the following we discuss other approaches that are related to the construction of almost-$C^1$ splines presented in this paper. Since almost-$C^1$ splines are smooth splines over quadrilateral meshes, we mostly focus the discussion on similar, quad-based spline constructions. Note that this is an extensive area of active research and we do not attempt to present an exhaustive overview; instead we point the readers to the comprehensive literature reviews recently presented in [25]. For instance, we focus here on constructions which build spline functions by smoothly joining polynomials on quadrilateral patches. This is in contrast to spline manifolds as in [19] where the functions are locally defined by composition of polynomials with suitably chosen blending functions, and which have been used to build bivariate [38, 70, 34] and trivariate splines [34] for IGA.

In the context of piecewise-polynomial splines over quadrilateral meshes, one can distinguish two types of smoothness – parametric and geometric. While parametrically smooth splines assume (at least locally) a joint parameter domain for neighboring elements, in which the smoothness is prescribed, geometric smoothness is defined directly in physical space. As we will see in the following, many constructions rely on both parametric and geometric smoothness. While geometric smoothness between Bézier or B-spline patches is a well-known concept in CAD [17, Chapter 8], it has only recently been used for simulations.

One possibility to create parametrically smooth surfaces from unstructured meshes is to use subdivision, see [16, 11, 58]. While subdivision is a process that can be defined entirely on a mesh, its limit surface can be interpreted as a piecewise polynomial spline. This limit surface has a finite representation in regular regions of the mesh (e.g. equivalent to bicubic B-splines in case of Catmull–Clark subdivision) and is composed of an infinite number of spline rings around extraordinary vertices, cf. [45]. When using subdivision surfaces in IGA, this peculiar property must be taken into account, e.g. when performing numerical quadrature [2]. Moreover, the approximation properties are in general sub-optimal, see [8, 1, 42, 66, 50, 69, 36].

It is also possible to employ parametric smoothness over a quadrilateral mesh composed of finitely many elements. While such a construction yields B-splines (or locally refined splines, such as T-splines) in structured regions of the mesh, singularities are introduced at extraordinary vertices. This phenomenon was studied and used in [47, 48] to create smooth but singular spline surfaces. In this setup, additional geometric smoothness has to be imposed at the extraordinary vertex. This results in smooth but degenerate Bézier patches, so-called D-patches. Such constructions were used for IGA in [43, 63, 10]. Even though the spline geometries are singular, the spaces possess many favorable properties for both design and analysis. They can be used to discretize high-order PDEs, see e.g. [71, 9], since they possess the required $H^2$-regularity properties, cf. [60]. Moreover, the spaces constructed in [63, 10] demonstrate optimal convergence under mesh refinement when applied to fourth-order problems.

Alternatively, one can increase the flexibility of smooth splines over quadrilateral meshes by creating polar singularities, which are the result of edges mapped onto single points. General $C^k$ smooth splines, for $k \geq 0$, were developed in [59] over singular Bézier patches and in [61] over general polar quad meshes. An explicit $C^1$ smooth construction was presented in detail in [57] and used for simulations on smooth, deforming surfaces in [65].

When constructing smooth splines over singular or polar configurations, additional geometric continuity has to be imposed at the singular or polar point to achieve the desired smoothness. Similar constructions can be employed, if parametric smoothness of higher order is enforced only for those edges that are away from extraordinary vertices, while at those edges near extraordinary vertices geometric smoothness is imposed, such as in [49, 52]. The dimension of such locally defined, geometrically smooth spline spaces over topological, mixed quadrilateral-triangle meshes was studied in [41]. Later, approximation properties could be shown for $G^1$ smooth isogeometric elements over planar, quadrilateral [31] and mixed [20] meshes.

In the following we give an overview of constructions that rely on geometric smoothness not in such a local, but in a global sense. Geometrically $C^1$ smooth isogeometric discretizations can be defined over bilinear Bézier elements [5], over bilinear spline patches [27], or over more general planar multi-patch domains [14, 28, 30]. Recently, constructions for $C^2$ smoothness over multi-patch domains were developed in [26, 32].

It has been shown in [14], that $C^1$ smooth spaces over multi-patch domains possess optimal approximation properties only in the case of so-called analysis-suitable $G^1$ parameterizations. While this condition can be fulfilled for planar domains following a reparameterization [29], for general planar multi-patch domains and surfaces it is not satisfied. Thus, to increase the geometric flexibility and allow for constructions over general domains, one may increase the polynomial degree locally [12, 13] or reduce the smoothness requirements by replacing exact $C^1$ smoothness by approximate $C^1$ smoothness [68]. In this paper we follow a similar approach, but instead of enforcing approximate $C^1$ smoothness along entire interfaces between patches, we impose approximate $C^1$ smoothness only at mesh edges near extraordinary vertices. In the following we give an overview of the most important properties of the almost-$C^1$ splines that we propose in this paper.

### 1.3. *Properties of almost-$C^1$ splines*

Given an unstructured mesh $\mathscr{T}$ consisting of quadrilaterals (we allow both planar and non-planar meshes of arbitrary topology), almost-$C^1$ splines are biquadratic splines on $\mathscr{T}$ that extend the construction developed in [64]. The following briefly outlines some features of our construction; they are elaborated upon later in the article.

- **Definable on fully unstructured meshes**: We allow all types of manifold quadrilateral meshes with no restrictions on the numbers or placements of extraordinary vertices (e.g., multiple extraordinary vertices per quadrilateral are allowed, boundary extraordinary vertices are allowed).
- **Well-conditioned B-spline-like basis**: The spline basis functions have several useful B-spline-like properties: partition of unity, non-negativity, local support and linear independence. The spline degree-of-freedom structure is simple and allows simple control of the geometry and functions at the boundary. Moreover, the splines are $C^1$ smooth (in an isogeometric sense) at all vertices of the mesh, and approximately smooth only across edges containing extraordinary vertices. Furthermore, the spline basis is non-singular and well-conditioned.
- **Easy implementation**: We utilize and extend the non-nested refinement process from [64] which is convergent and nested, when restricted to the boundaries. As outlined in [64], this has two benefits which improve upon nested refinements of spline spaces of mixed smoothness (e.g., see [63, 67]). The non-nestedness of the refinement allows us to "shrink" the neighbourhoods of approximate $C^1$ smoothness and leads to a very simplified computer implementation and, furthermore, in the limit of infinite refinements, converges to a smooth limit surface. At the same time, the refinement process leaves the spline invariant on the mesh boundary; this is especially useful if the boundaries of a spline geometry are composed of special curves such as conic sections.

  Numerical tests indicate that the spline spaces also demonstrate optimal approximation behaviour in the $L^2$, $H^1$ and $H^2$ norms for second- and fourth-order problems under mesh refinement. Conceptually, this approach can be seen as an amalgamation of the "design" and "analysis" philosophies from [63] – it offers the ease of working with the design space while also being suitable for analysis.
- **Lowest-order unstructured spline construction for fourth-order problems**: Our spline construction only uses biquadratic polynomial pieces and is thus the lowest-order unstructured spline construction suitable for fourth-order problems. Note that it is well-known that for certain problems (e.g., from structural mechanics) higher polynomial degrees might help alleviate *locked* or non-convergent approximations; similar ad-hoc solutions can also be derived for lower polynomial degrees (for instance, by utilizing reduced quadrature).

### 1.4. *Outline of the paper*

In Section 2 we introduce the relevant notation for the unstructured quadrilateral meshes that we focus on. Then, in Section 3 we present the construction of unstructured biquadratic splines over such meshes, culminating in the definition of almost-$C^1$ splines in Section 3.5. Their useful properties are presented in Proposition 3.8. Next, in Section 4, we discuss the (non-nested) refinement of almost-$C^1$ spline spaces and geometries. Proposition 4.3 characterizes our refinement rules. Finally, in Section 5 we present some numerical tests focusing on the analysis-suitability of the proposed B-splines, including the Scordelis–Lo thin shell benchmark, a Cahn–Hilliard problem on a closed surface and the analysis of a Laplace–Beltrami eigenvalue problem.

## 2. Unstructured quadrilateral meshes

We are interested in solving scalar and vector-valued PDEs on complex, 2-dimensional geometries of arbitrary topology. For instance, planar geometries in $\mathbb{R}^2$ or surfaces in $\mathbb{R}^3$. Then, splines defined on unstructured meshes can help create such complex geometries, and can thereafter be used for numerically solving PDEs on them. We focus here on unstructured quadrilateral meshes. In this section, we define some relevant notation for such meshes, which are the basis for defining almost-$C^1$ splines; see Figure 2(a) for reference.

Before we begin, we would like to point out that we consider $\mathscr{T}$ as a topological construct only — in general, the quadrilaterals in $\mathscr{T}$ will not be assumed to occupy a common parametric domain. Similarly, the meshes are not restricted to be planar or of trivial topology either. Representations as in Figure 2(a) will only be for the purpose of specifying the connectivity of the different quadrilaterals with each other. As such, spline geometries and spline functions on those geometries will be built by appropriately selecting the degrees of freedom for splines on $\mathscr{T}$. An example corresponding to the mesh in Figure 2(a) is shown in Figure 2(c); see Section 3.5 for details on the construction of almost-$C^1$ splines. Note that the spline construction presented there is based on local geometric data around extraordinary vertices. However, the construction can easily be generalized to a purely topological one, as explained in Remark 3.4 and Appendix B.

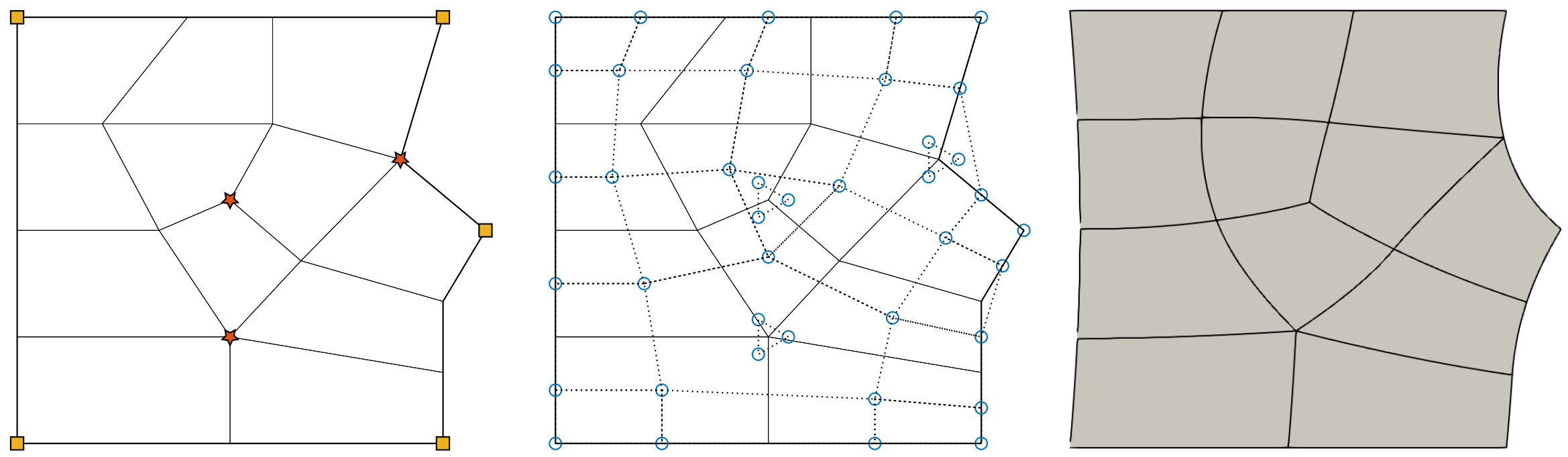

(a) mesh topology (b) spline dof structure (c) example geometry

Figure 2: An example of almost-$C^1$ splines as defined in Section 3, with the underlying mesh in (a), the dof-structure in (b) and an example geometry in (c).

We will denote all *quadrilateral meshes* with $\mathscr{T}$. We assume that $\mathscr{T}$ is without any hanging nodes and that the interiors of all quadrilaterals are disjoint. Each quadrilateral in $\mathscr{T}$ will be called a *face of* $\mathscr{T}$, or simply a *face*. More generally, for $k \leq 2$, the $k$-dimensional geometric components of the mesh will be collected in sets $\mathscr{T}_k$. That is, vertices in $\mathscr{T}_0$, edges in $\mathscr{T}_1$, and faces in $\mathscr{T}_2$.

We assume that $\mathscr{T}$ is such that each edge is shared by at most two faces of the mesh. If any edge is contained in only one face then it is called a *boundary edge*, otherwise it is called an *interior edge*; boundary edges have been displayed with slightly thicker lines in Figure 2(a). Any vertices that lie on a boundary edge are called *boundary vertices*, otherwise they are called *interior vertices*. We denote the sets of boundary edges and vertices with $\mathscr{T}_1^B$ and $\mathscr{T}_0^B$, respectively. The set of interior vertices and edges will be denoted by $\mathring{\mathscr{T}}_k$ for $k = 0$ and 1, respectively. We also assume that there are no 'kissing vertices' in the mesh. That is, for any two faces $\sigma, \sigma'$ that share a common vertex $\gamma$, there is a sequence of faces $\sigma_0, \dots, \sigma_\ell$ that all contain $\gamma$ such that $\sigma_0 = \sigma$, $\sigma_i = \sigma'$, and $\sigma_i \cap \sigma_{i-1} \in \mathring{\mathscr{T}}_1$ for $i = 1, \dots, \ell$.

We always assume faces and edges to be closed sets. Thus, we define the *valence* of a vertex, edge or face of $\mathscr{T}$ to be the number of faces that contain it. For $\phi \in \mathscr{T}_k$, $k \leq 2$, we will denote the valence of $\phi$ with $\mu_\phi$. In particular,

- the valence $\mu_\sigma$ of any face $\sigma$ is exactly 1 since each face contains itself;
- the valence $\mu_\tau$ of a boundary or interior edge $\tau$ is 1 or 2, respectively, by the above definitions.

Vertices $\gamma$ of $\mathscr{T}$ will be called *extraordinary vertices* if they are interior vertices with valences $\mu_\gamma \neq 4$, or if they are boundary vertices with valences $\mu_\gamma > 2$. Vertices of valence 1 will be called *corner vertices* and will be collected in the set $\mathscr{T}_0^C \subseteq \mathscr{T}_0^B = \mathscr{T}_0 \backslash \mathring{\mathscr{T}}_0$. For visual consistency, we will always denote the extraordinary and corner vertices of a mesh by respectively placing red stars and yellow squares on them; see Figure 2(a). Moreover, we call all edges that contain one or more extraordinary vertices as *spoke edges* and all faces that contain one or more extraordinary vertices as *extraordinary faces*. Otherwise, we call them regular edges and regular faces, respectively. We denote the set of all extraordinary vertices and faces with $\mathscr{T}_0^E$ and $\mathscr{T}_2^E$, respectively.

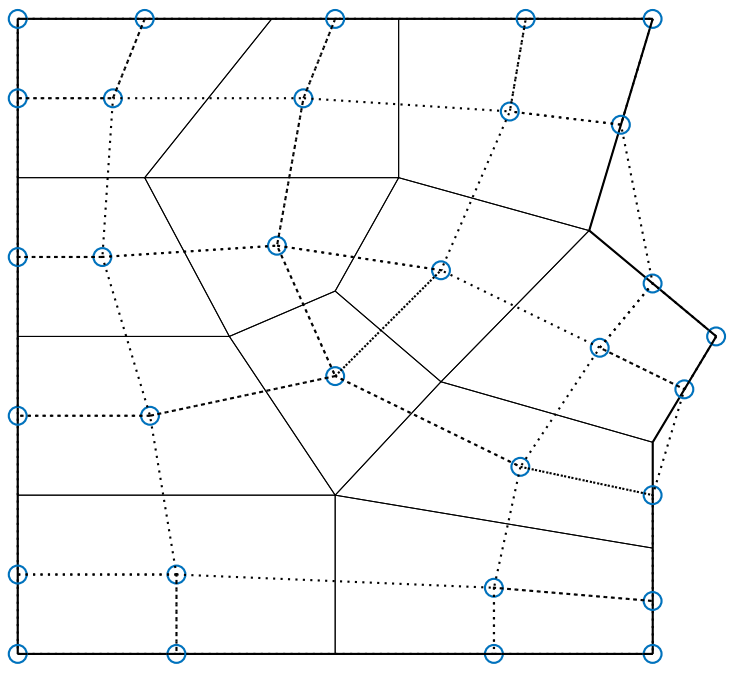

Figure 3: The degree of freedom structure for the spline space $\mathscr{B}^*$ on the mesh $\mathscr{T}$ from Figure 2(a). There is one dof associated to each face, boundary edge and corner vertex of the mesh $\mathscr{T}$.

## 3. Unstructured biquadratic splines

In this section we present the construction of splines over unstructured quadrilateral meshes $\mathscr{T}$ as defined above. The construction is based on two steps. First, we summarize a biquadratic spline basis construction from [64] which depends purely on topological properties of the mesh (Sections 3.1–3.3). These functions span the space $\mathscr{B}^*$, which is biquadratic on every face, $C^0$ at all extraordinary points and across all spoke edges and $C^1$ across all other edges of the mesh. In the second step the basis functions of $\mathscr{B}^*$, which have support on extraordinary faces, are modified such that they have vanishing value and gradient at each extraordinary point. Consequently, three additional functions for each extraordinary vertex are introduced that locally span all linear functions. These functions can be defined using some geometric information of the mesh. The resulting functions constitute the almost-$C^1$ splines over the mesh $\mathscr{T}$, spanning the space $\mathscr{B}$ (see Section 3.5). The almost-$C^1$ splines are then $C^1$ smooth at all vertices and composed of modified functions from $\mathscr{B}^*$ as well as three new functions associated to each extraordinary vertex.

### *3.1. Degree-of-freedom structure for $\mathscr{B}^*$*

The degrees of freedom, or *dofs* in short, corresponding to $\mathscr{B}^*$ are divided into three categories; the spline construction will be specified for each category separately.

- Face dofs: We associate one degree of freedom to each face of $\mathscr{T}$, i.e., one degree of freedom for each member of $\mathscr{T}_2$.
- Boundary edge dofs: We associate one degree of freedom to each boundary edge of $\mathscr{T}$, i.e., one degree of freedom for each member of $\mathscr{T}_1^B$.
- Corner vertex dofs: We associate one degree of freedom to each corner vertex of $\mathscr{T}$, i.e., one degree of freedom for each member of $\mathscr{T}_0^C$.

We create a basis for $\mathscr{B}^*$ by associating one function to each dof, which we summarize using the index set

$$\mathscr{I}^* := \mathscr{T}_2 \cup \mathscr{T}_1^B \cup \mathscr{T}_0^C \,. \tag{1}$$

For visual consistency, in the topological description (as in Figure 2(b)), we will denote each dof by placing an unfilled blue circle on the associated face/boundary edge/corner vertex of $\mathscr{T}$, and the connectivity of the dofs will be denoted with thin dotted lines; see Figure 3.

Following this categorization of the dofs, we will also call the associated splines *face*, *boundary edge* and *corner vertex splines*, respectively. Denote these spline functions as $B^*_\phi$, $\phi \in \mathscr{I}^*$. The associated spline space on $\mathscr{T}$ is then going to be defined as

$$\mathscr{B}^* := \operatorname{span}\left(B^*_\phi \ : \ \phi \in \mathscr{I}^*\right) ,$$

where $\mathscr{I}^*$ is defined as in (1).

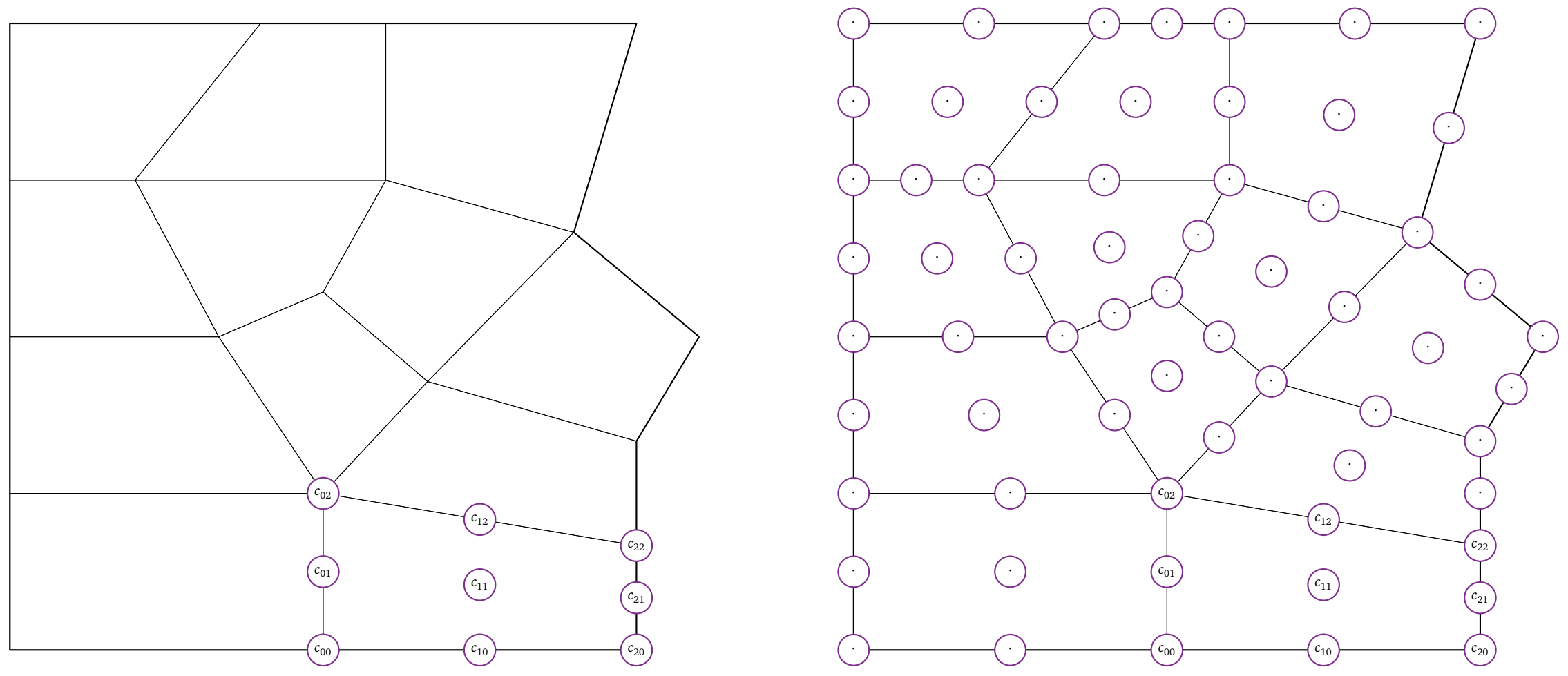

(a) element-local coefficients

(b) extension to neighbours by 0

Figure 4: Here we take a submesh of the one from Figure 2(a) to illustrate how the extraction coefficients will be graphically denoted. Given a face $\sigma \in \mathcal{T}$ and $\phi$ such that $\phi \cap \sigma = \phi$, $B^*_\phi$ will be defined graphically by specifying its Bernstein–Bézier coefficients on $\sigma$; see figure (a). This description will be extended to neighbouring faces of $\mathcal{T}$ with zero coefficients denoted here with a $\cdot$ in figure (b). See Section 3.2 for further elaboration.

In Section 3.3, we define the B-splines $B^*_\phi$, $\phi \in \mathscr{I}^*$. This will be done by specifying their local polynomial descriptions in terms of Bernstein–Bézier polynomials on quadrilaterals via so-called *extraction matrices*. Extraction matrices specify how the face-local Bernstein–Bézier polynomials can be linearly combined to yield the face-local description of a spline basis function; these were introduced in [7, 51] and have been used, for instance, for defining splines on unstructured quadrilateral meshes in [61, 63, 67] and for multi-degree splines in [62, 55]. We explain our extraction matrix convention in Section 3.2, and the B-spline definitions are subsequently presented in Section 3.3.

### 3.2. Extraction matrix convention

We will graphically denote the extraction matrix for each B-spline. For all B-splines, this will be done by specifying their Bernstein–Bézier coefficients on a single face of $\mathcal{T}$; the face-local description will be extended to the neighbours. We elaborate upon this convention here and use Figure 4 for reference. For visual consistency and to differentiate local Bernstein–Bézier coefficients from spline dofs, all local polynomial coefficients will be displayed inside magenta coloured disks.

Consider a face $\sigma \in \mathcal{T}_2$ and let $B^*_\phi$ be a spline associated either to $\sigma$ or to a boundary edge/corner vertex that belongs to $\sigma$, i.e., $\phi \cap \sigma = \phi$. Then, we will present the definition of $B^*_\phi$ graphically as on the left in Figure 4(a) by specifying 9 coefficients $\{c_{jk}[B^*_\phi;\sigma]\}_{j,k=0}^2$. The shown coefficients are to be interpreted as defining the following local spline description,

$$B^*_\phi\big|_\sigma = \sum_{j,k=0}^{2} c_{jk}[B^*_\phi;\sigma]\, b^0_{jk,\square}\,, \tag{2}$$

where $b^0_{jk,\square}$ is the $(j,k)$-th biquadratic Bernstein–Bézier polynomial defined on $\sigma$ by interpreting it as the unit square $[0,1]^2$ (the origin is placed at the corner with the coefficient $c_{00}$) and local coordinates $\xi := (u,v)$,

$$b^0_{jk,\square}(\xi) := \binom{2}{j}\binom{2}{k}(1-u)^j(1-v)^k u^j v^k\,. \tag{3}$$

Thus, $c_{jk}[f;\sigma]$ denotes the $(j,k)$-th coefficient of the function $f$ restricted to the face $\sigma$.

Finally, the face-local descriptions specified as above on any face $\sigma$ are extended to any neighbouring face $\sigma'$ with the help of zero coefficients; see Figure 4(b). Note that on a planar quadrilateral mesh such a representation

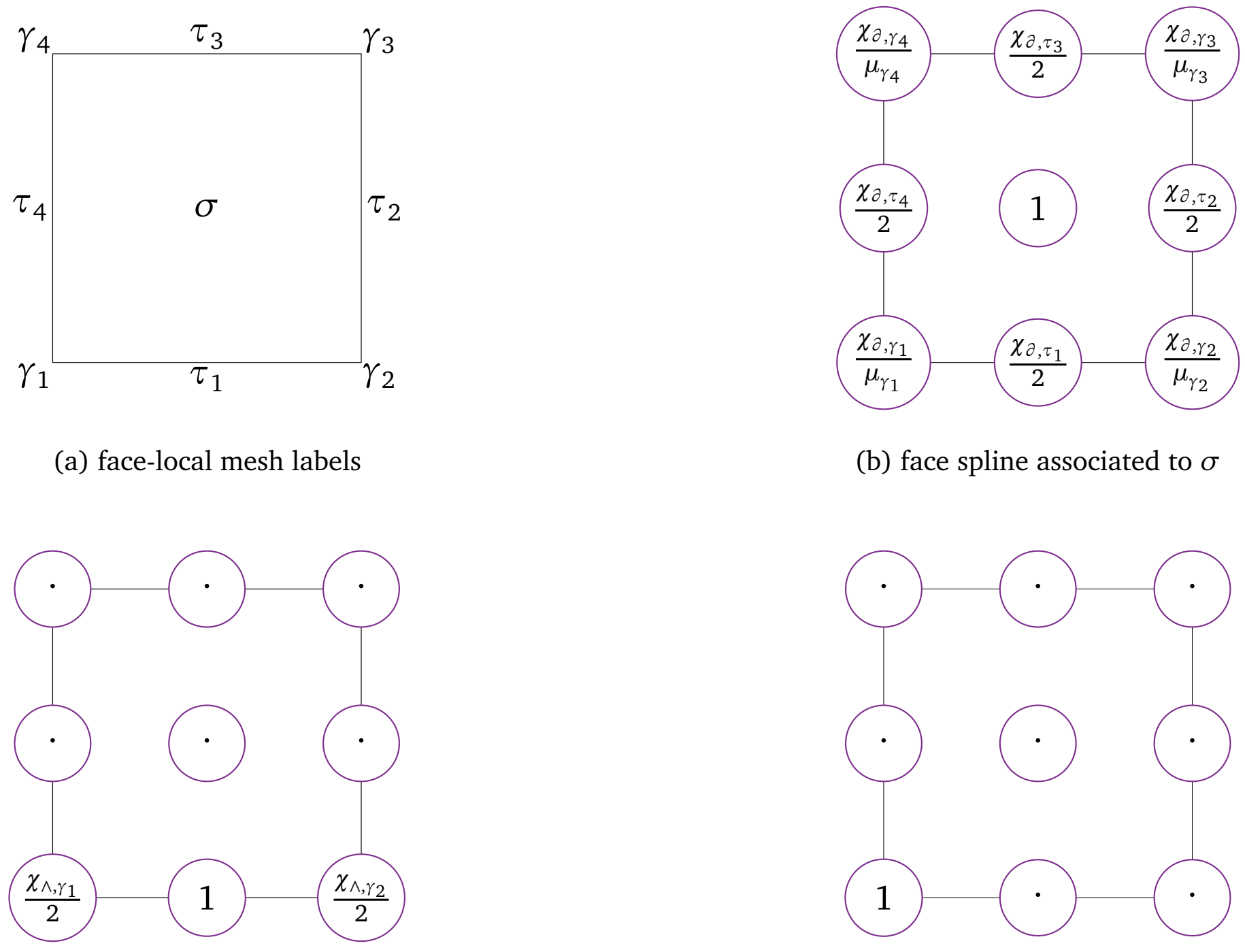

(a) face-local mesh labels

(b) face spline associated to $\sigma$

(c) boundary-edge spline associated to $\tau_1$

(d) corner-vertex spline associated to $\gamma_1$

Figure 5: This figure defines all B-splines that may be associated to a face, boundary edge or a corner vertex. In (a), the local mesh neighbourhood of a face of $\mathscr{T}$ are shown with vertices $\gamma_i$ and edges $\tau_i$. In (b), the face-local coefficients of a B-spline associated to $\sigma$ are presented. In (c), assuming that $\tau_1$ is a boundary edge, the face-local coefficients of a B-spline associated to $\tau_1$ are presented. Finally, in (d), assuming that $\gamma_1$ is a corner vertex, the face-local coefficients of a B-spline associated to $\gamma_1$ are presented.

corresponds to representing the dof-structure of the $C^0$-smooth piecewise biquadratic space using domain points, as e.g. in [35]. With this graphical convention in place, let us now define the three different types of splines $B^*_\phi$, $\phi \in \mathscr{I}^*$, in the following section.

### *3.3. The element-local spline representations*

**Face splines**

Consider a face $\sigma \in \mathscr{T}_2$ and let $B^*_\sigma$ be the B-spline associated to its dof. Let the edges and vertices of $\sigma$ be numbered as shown at the top in Figure 5(a). Then, the face-local coefficients for $B^*_\sigma$ are defined as in Figure 5(b). The figure uses a *boundary flag* for vertices and edges which is a characteristic function defined as follows for $\phi \in \mathscr{T}_k$, $k = 0, 1$,

$$\chi_{\partial,\phi} = \begin{cases} 1\,, & \phi \in \mathring{\mathscr{T}}_k\,, \\ 0\,, & \text{otherwise}\,. \end{cases} \tag{4}$$

**Boundary edge splines**

Consider a face $\sigma \in \mathscr{T}_2$. Let the edges and vertices of $\sigma$ be numbered as shown at the top in Figure 5(a), and let $\tau_1$ be a boundary edge and $B^*_{\tau_1}$ the B-spline associated to its dof. Then, the face-local coefficients for $B^*_{\tau_1}$ are defined as in Figure 5(c). The figure uses a *corner flag* for vertices which is a characteristic function defined as follows for $\gamma \in \mathscr{T}_0$,

$$\chi_{\wedge,\gamma} = \begin{cases} 0\,, & \gamma \in \mathscr{T}_0^C\,, \\ 1\,, & \text{otherwise}\,. \end{cases} \tag{5}$$

**Corner vertex splines**

Consider a face $\sigma \in \mathscr{T}_2$. Let the edges and vertices of $\sigma$ be numbered as shown at the top in Figure 5(a). Let $\gamma_1$ be a boundary vertex that has been chosen to be a corner vertex, and let $B^*_{\gamma_1}$ the B-spline associated to its dof. Then, the face-local coefficients for $B^*_{\gamma_1}$ are defined as in Figure 5(d).

*3.4. Spline functions*

Given the above B-spline definitions, we can create a spline function $f^* \in \mathscr{B}^*$ by linearly combining the B-splines with coefficients $f^*_\phi \in \mathbb{R}$,

$$f^* := \sum_{\phi \in \mathscr{I}^*} f^*_\phi B^*_\phi \,. \tag{6}$$

Following the piecewise definition of B-splines, the above definition is also interpreted in a piecewise manner. That is, for any $\sigma \in \mathscr{T}_2$ and local coordinates $\boldsymbol{\xi}$ as in Equation (3),

$$f^*\big|_\sigma(\boldsymbol{\xi}) = \sum_{\phi \in \mathscr{I}^*} f^*_\phi B^*_\phi\big|_\sigma(\boldsymbol{\xi}) \,, \tag{7}$$

where $B^*_\phi\big|_\sigma$ follow from Section 3.3. Similar to (2), we thus have

$$f^*\big|_\sigma(\boldsymbol{\xi}) = \sum_{j,k=0}^{2} c_{jk}[f^*;\sigma] b_{jk,\Box}(\boldsymbol{\xi}) \,,$$

where

$$c_{jk}[f^*;\sigma] = \sum_{\phi \in \mathscr{I}^*} f^*_\phi c_{jk}[B^*_\phi;\sigma] \,.$$

These kinds of linear combinations can be used, for instance, to create a bivariate spline geometry $\boldsymbol{x}^* \in (\mathscr{B}^*)^d$. This can be done by choosing appropriate control points $\boldsymbol{x}^*_\phi \in \mathbb{R}^d$, $\phi \in \mathscr{I}^*$, and defining

$$\boldsymbol{x}^* := \sum_{\phi \in \mathscr{I}^*} \boldsymbol{x}^*_\phi B^*_\phi \,. \tag{8}$$

We also have a local representation for the geometry $\boldsymbol{x}^*$,

$$\boldsymbol{c}_{jk}[\boldsymbol{x}^*;\sigma] = \sum_{\phi \in \mathscr{I}^*} \boldsymbol{x}^*_\phi c_{jk}[B^*_\phi;\sigma].$$

The functions $\{B^*_\phi\}_{\phi \in \mathscr{I}^*}$ are linearly independent. They also form a non-negative, local partition of unity, thus the coefficients $\boldsymbol{x}^*_\phi$ can be seen as classical spline control points. See [64] for a discussion of other properties of $\mathscr{B}^*$.

*3.5. The almost-$C^1$ splines $\mathscr{B}$*

In the following we discuss how $\mathscr{B}^*$ and its basis can be modified to build almost-$C^1$ splines spanning the space $\mathscr{B}$, which is of interest for higher order problems. We focus on one specific construction which uses geometric data, i.e., which is based on an underlying geometry mapping $\boldsymbol{x}^* \in (\mathscr{B}^*)^d$, or equivalently on underlying control points $\boldsymbol{x}^*_\phi \in \mathbb{R}^d$, $\phi \in \mathscr{I}^*$, as well as on extraordinary-vertex normals $\boldsymbol{n}_\gamma$, for each extraordinary vertex $\gamma \in \mathscr{T}_0^E$. Note that the construction depends only on a small neighborhood of each extraordinary vertex, that is, on the control points corresponding to extraordinary faces $\phi \in \mathscr{T}_2^E \subset \mathscr{I}^*$. A completely geometry independent construction is also possible, as developed in Appendix B.

First, the set of degrees of freedom $\mathscr{I}$ for $\mathscr{B}$ is given by the dofs from $\mathscr{I}^*$ enriched by three additional dofs for each extraordinary vertex $\gamma \in \mathscr{T}_0^E$. More precisely, we set

$$\begin{aligned} \mathscr{I} = &\left\{\sigma \in \mathscr{T}_2 \;:\; \exists \gamma \in \mathscr{T}_0 \setminus \mathscr{T}_0^E \text{ s.t. } \gamma \in \sigma\right\} \\ &\bigcup \left\{\tau \in \mathscr{T}_1^B \;:\; \exists \gamma \in \mathscr{T}_0 \setminus \mathscr{T}_0^E \text{ s.t. } \gamma \in \tau\right\} \\ &\bigcup \mathscr{T}_0^C \bigcup \left(\mathscr{T}_0^E \times \{1,2,3\}\right), \end{aligned} \tag{9}$$

and describe a basis $B_\phi$, $\phi \in \mathscr{I}$, in the following. Note that there are no splines in $\mathscr{B}$ associated to those faces and boundary edges where all incident vertices are extraordinary. The space $\mathscr{B}$ is defined as the span of $B_\phi$,

$$\mathscr{B} := \operatorname{span}\left(B_\phi \;:\; \phi \in \mathscr{I}\right) \,.$$

**B-splines on regular faces are identical**

Let $\phi \in \mathscr{I}^*$ be a dof of the spline space $\mathscr{B}^*$. Then we define the corresponding basis function $B_\phi$ of $\mathscr{B}$ such that

$$B_\phi|_\sigma := B^*_\phi|_\sigma\,, \quad \forall \sigma \in \mathscr{T}_2 \backslash \mathscr{T}_2^E\,.$$

**B-splines on extraordinary faces are subdivided and truncated**

For $\phi \in \mathscr{I}^*$ and $\sigma \in \mathscr{T}_2^E$, we will define $B_\phi|_\sigma$ by modifying (i.e., subdividing and truncating) the local representations $B^*_\phi|_\sigma$; the modification will impose vanishing values and derivatives for $B_\phi$ at all $\gamma \in \mathscr{T}_0^E$. If the modifications imply that $B_\phi \equiv 0$, then $\phi$ will not be a dof for the spline space $\mathscr{B}$; otherwise, $B_\phi$ will be the basis function corresponding to the dof $\phi$. In particular, the former situation will arise only if all vertices incident upon $\phi$ are extraordinary.

First, for any $\sigma \in \mathscr{T}_2^E$, consider the local representation of $B^*_\phi|_\sigma$,

$$B^*_\phi|_\sigma = \sum_{j,k=0}^{2} c_{jk}[B^*_\phi;\sigma]b^0_{jk,\square}\,,$$

and define the matrix of coefficients $c[f;\sigma] := \left(c_{jk}[f;\sigma]\right)_{jk}$. Let the four vertices of $\sigma$ be ordered 1 through 4 in counter-clockwise manner starting from the bottom-left, as in Figure 5(a). Then, we locally define $B_\phi$ to be a modified representation of $B^*_\phi$ as

$$B_\phi|_\sigma := \sum_{j,k=0}^{3} \hat{c}_{jk}[B_\phi;\sigma]b^1_{jk,\square}\,, \tag{10}$$

where $b^1_{jk,\square}$ is the $(j,k)$-th $C^1$ tensor-product B-spline corresponding to the knot vector $(0,0,0,\frac{1}{2},1,1,1)$ in both parametric directions, and

$$\hat{c}[B_\phi;\sigma] := \left(\mathrm{K}c[B^*_\phi;\sigma]\mathrm{K}^T\right)\bigodot_{i=1}^{4}\mathrm{T}_i\,,$$

with $\odot$ signifying the Hadamard product (or, element-wise product) of matrices. In the above, K is the univariate B-spline knot-insertion matrix that takes a quadratic Bézier element $[0,1]$ and inserts a single knot at 0.5,

$$\mathrm{K} = \begin{bmatrix} 1 & & \\ \frac{1}{2} & \frac{1}{2} & \\ & \frac{1}{2} & \frac{1}{2} \\ & & 1 \end{bmatrix}\,, \tag{11}$$

and $\mathrm{T}_i$ is the truncation matrix associated to the $i$-th vertex of $\sigma$ as ordered above. If the $i$-th vertex is not extraordinary, $\mathrm{T}_i$ is defined to be a $4\times 4$ matrix with all entries equal to 1. Else, if the $i$-th vertex is extraordinary, $i = 1,\ldots,4$, we respectively define

$$\mathrm{T}_1 = \begin{bmatrix} 0 & 0 & 1 & 1 \\ 0 & 0 & 1 & 1 \\ 1 & 1 & 1 & 1 \\ 1 & 1 & 1 & 1 \end{bmatrix}\,,\quad \mathrm{T}_2 = \begin{bmatrix} 1 & 1 & 1 & 1 \\ 1 & 1 & 1 & 1 \\ 0 & 0 & 1 & 1 \\ 0 & 0 & 1 & 1 \end{bmatrix}\,,\quad \mathrm{T}_3 = \begin{bmatrix} 1 & 1 & 1 & 1 \\ 1 & 1 & 1 & 1 \\ 1 & 1 & 0 & 0 \\ 1 & 1 & 0 & 0 \end{bmatrix}\,,\quad \mathrm{T}_4 = \begin{bmatrix} 1 & 1 & 0 & 0 \\ 1 & 1 & 0 & 0 \\ 1 & 1 & 1 & 1 \\ 1 & 1 & 1 & 1 \end{bmatrix}\,. \tag{12}$$

Note that if all vertices of $\sigma$ are extraordinary, $\hat{c}[B_\phi;\sigma] = \mathbf{0}$ since $\odot_{i=1}^4 \mathrm{T}_i = \mathbf{0}$. In particular, if $\hat{c}[B_\phi;\sigma] = \mathbf{0}$ for $\sigma$ such that $\sigma \cap \phi = \phi$, then $B_\phi$ will be globally zero since the corresponding $B^*_\phi$ is supported only on extraordinary faces and is entirely truncated on each face in its support. Furthermore, for all $\phi \in \mathscr{I}^*$ such that $\phi \cap \gamma = \emptyset$ for all $\gamma \in \mathscr{T}_0^E$, we have $B_\phi = B^*_\phi$. That is, all regular (face, boundary edge and corner) B-splines are unchanged by the subdivision and truncation.

*Remark* 3.1. A similar construction can be achieved if the functions $b^1_{jk,\square}$ in (10) are replaced by bicubic polynomials and if K is replaced by a degree elevation matrix instead of the knot insertion matrix, leading to almost-$C^1$ splines being bicubic polynomials on extraordinary faces. Details can be found in Appendix A.

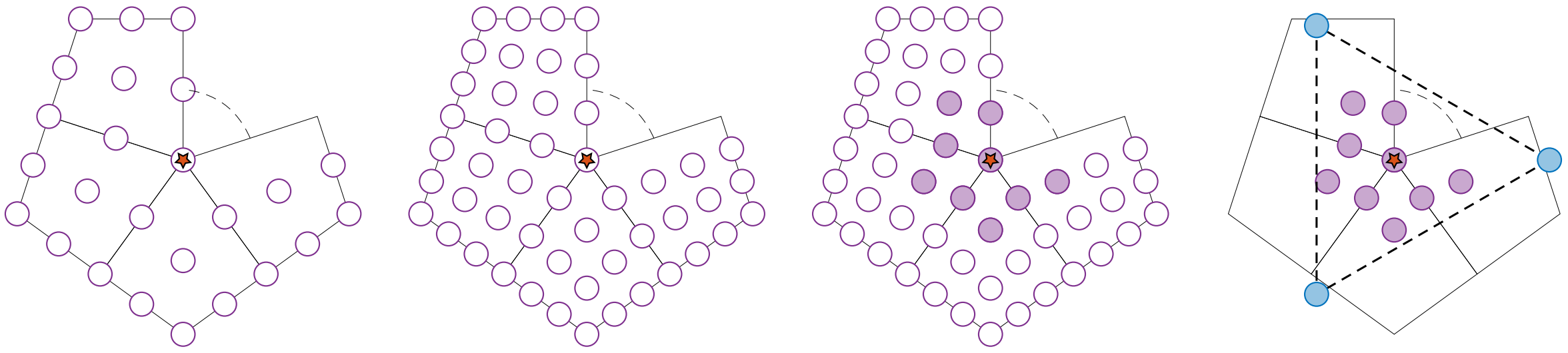

(a) coefficients $c_{jk}[\boldsymbol{x}^*;\sigma_i]$ (b) coefficients $\hat{c}_{jk}[\boldsymbol{x}^*;\sigma_i]$ (c) projected coefficients (d) control triangle

Figure 6: An extraordinary vertex $\gamma \in \mathscr{T}_0^E$ surrounded by its 1-ring, composed of the extraordinary faces $\sigma_1, \sigma_2, \ldots, \sigma_\mu$. The coefficients $c_{jk}[\boldsymbol{x}^*;\sigma_i]$ as well as the once refined coefficients $\hat{c}_{jk}[\boldsymbol{x}^*;\sigma_i]$ of $\boldsymbol{x}^*$ within the 1-ring neighbourhood of $\gamma$ are depicted in (a) and (b), respectively. In (c) the relevant coefficients needed to define $B_{(\gamma,\nu)}$, as in (13), are highlighted in purple. The control triangle covering the relevant coefficients is visualized in (d).

**Extraordinary vertex splines $B_{(\gamma,i)}$ are added**

Let $\gamma \in \mathscr{T}_0^E$ be an extraordinary vertex of valence $\mu$ and let the faces around it be denoted by $\sigma_1, \sigma_2, \ldots, \sigma_\mu$. Given a prescribed normal direction $\boldsymbol{n}_\gamma$ at the vertex, cf. Remark 3.3, we introduce a tangent plane through $\boldsymbol{x}^*|_\gamma$ which is orthogonal to the vector $\boldsymbol{n}_\gamma$. Then, the 1-ring $\boldsymbol{x}^*|_{\bigcup_{i=1}^{\mu}\sigma_i}$ is projected orthogonally onto the tangent plane. We denote this projection by $P_\gamma$. Since the construction is affine invariant, this is equivalent to projecting the Bernstein–Bézier coefficients of the function $\boldsymbol{x}^*|_{\sigma_i}$ for each face $\sigma_i$, i.e.,

$$P_\gamma(\boldsymbol{x}^*)|_{\sigma_j} = \sum_{j,k=0}^{3} P_\gamma(\hat{\boldsymbol{c}}_{jk}[\boldsymbol{x}^*;\sigma_j])b^1_{jk,\square}.$$

*Remark* 3.2. Note that $P_\gamma(\boldsymbol{x}^*)|_{\sigma_j}$ could have self-intersections but that is irrelevant for our construction where we only need the projections of 4 specific element-local coefficients in the vicinity of $\gamma$ to be regular; see Figure 6, Equation 13 and Remark 3.3.

Thus, $P_\gamma$ results in a planar configuration of faces around the extraordinary vertex, see Figure 6. We select the local control points

$$\boldsymbol{c}^i_{jk} := P_\gamma(\hat{\boldsymbol{c}}_{jk}[\boldsymbol{x}^*;\sigma_i]), \quad \text{with } i \in \{1,\ldots,\mu\},\ (j,k) \in \{0,1\}^2\,, \tag{13}$$

which are relevant for the basis construction. These local control points $\boldsymbol{c}^i_{jk}$, as highlighted in Figure 6(c), are then covered by a triangle $(\boldsymbol{x}_{(\gamma,1)}, \boldsymbol{x}_{(\gamma,2)}, \boldsymbol{x}_{(\gamma,3)})$, shown in Figure 6(d). The triangle is selected such that it is the smallest triangle that contains all points $\boldsymbol{c}^i_{jk}$ in its interior (or of similar size to the smallest triangle, cf. Remark 3.5). We denote the barycentric coordinates of a point $\boldsymbol{c}$ with respect to the triangle $(\boldsymbol{x}_{(\gamma,1)}, \boldsymbol{x}_{(\gamma,2)}, \boldsymbol{x}_{(\gamma,3)})$ by $\lambda_1(\boldsymbol{c})$, $\lambda_2(\boldsymbol{c})$ and $\lambda_3(\boldsymbol{c})$, i.e., we have

$$\boldsymbol{c} = \lambda_1(\boldsymbol{c})\boldsymbol{x}_{(\gamma,1)} + \lambda_2(\boldsymbol{c})\boldsymbol{x}_{(\gamma,2)} + \lambda_3(\boldsymbol{c})\boldsymbol{x}_{(\gamma,3)}\,,$$

with $\lambda_1(\boldsymbol{c}) + \lambda_2(\boldsymbol{c}) + \lambda_3(\boldsymbol{c}) = 1$. By construction, we have $\lambda_\nu(\boldsymbol{c}^i_{jk}) \geq 0$ for all $\nu \in \{1,2,3\}$, $(j,k) \in \{0,1\}^2$ and $i \in \{1,\ldots,\mu\}$.

The triangle $(\boldsymbol{x}_{(\gamma,1)}, \boldsymbol{x}_{(\gamma,2)}, \boldsymbol{x}_{(\gamma,3)})$ serves as a control triangle to determine the coefficients of the three new basis functions $B_{(\gamma,\nu)}$, with $\nu \in \{1,2,3\}$. The coefficients of $B_{(\gamma,\nu)}|_{\sigma_i}$ are given as the barycentric coordinates corresponding to the local control point $\boldsymbol{x}_{(\gamma,\nu)}$, i.e., the function $B_{(\gamma,\nu)}$ is defined to be

$$B_{(\gamma,\nu)}|_{\sigma_i} = \sum_{j,k=0}^{3} \hat{c}_{jk}[B_{(\gamma,\nu)};\sigma_i]b^1_{jk,\square}\,,$$

with

$$\hat{c}_{jk}[B_{(\gamma,\nu)};\sigma_i] = \begin{cases} \lambda_\nu(\boldsymbol{c}^i_{jk})\,, & (j,k) \in \{0,1\}^2\,, \\ 0\,, & \text{otherwise}\,. \end{cases}$$

Moreover, $B_{(\gamma,\nu)}|_\sigma = 0$ for all $\sigma \notin \{\sigma_1, \sigma_2, \ldots, \sigma_\mu\}$.

Once the above process is repeated for all extraordinary points, the geometry description is updated from $\boldsymbol{x}^*$ to $\boldsymbol{x} \in \mathscr{B}^d$, with the latter defined to be

$$\boldsymbol{x} = \sum_{\phi \in \mathscr{I}} \boldsymbol{x}_\phi B_\phi \; . \tag{14}$$

*Remark* 3.3. The prescribed normal vector $\boldsymbol{n}_\gamma$ must be given such that the projection is well-defined in a neighbourhood of the extraordinary vertex, i.e., $\boldsymbol{n}_\gamma \neq \mathbf{0}$, and the projected surface $P_\gamma(\boldsymbol{x}^*)$ is regular *in a neighbourhood* of $\gamma$. This implies a mild regularity assumption on the geometry $\boldsymbol{x}^*$. If no normal vector is prescribed, a suitable normal vector can be constructed e.g. through averaging of local normals in a neighbourhood of $\gamma$.

*Remark* 3.4. The projection onto the tangent plane prescribed by $\boldsymbol{n}_\gamma$ introduces a dependence of the almost-$C^1$ splines on the geometry $\boldsymbol{x}^*$. This is not necessary. Instead, for all valences that appear in the mesh $\mathscr{T}$, one can also prescribe regular templates that enforce coplanarity in a geometry independent manner. The resulting coefficients for such templates are given in Appendix B for some common valences.

*Remark* 3.5. The choice of the control triangle has no effect on the resulting space but only on the properties of the constructed basis. Thus, the properties desired from the basis can also inform the choice of the control triangle. We briefly mention two cases here: one concerning the conditioning of the basis, and another concerning basis functions at the boundary. First, regarding conditioning, the control triangle should be chosen to be sufficiently regular (such that the smallest angle is bounded away from zero) and not too large. See [56] where control triangles were used to specify values and derivatives for splines over triangulations. Second, in case of an extraordinary vertex $\gamma$ at the boundary, one can select the control triangle to obtain a basis that behaves like the usual $C^1$ quadratic B-spline basis when restricted to the boundary. Let $\tau_1, \tau_2 \in \mathscr{T}_1^B$ be the two spoke edges containing $\gamma$, and let the projected control points that correspond to one or both of these boundary edges be $\mathbf{c}_1^1$, $\mathbf{c}_1^2$ and $\mathbf{c}_0 := \frac{\mathbf{c}_1^1+\mathbf{c}_1^2}{2}$. Then all other projected control points in (13) lie on one side of the line through $\mathbf{c}_1^1$ and $\mathbf{c}_1^2$. Thus, choosing a control triangle such that this line contains one of its edges, only two of the basis functions $B_{(\gamma,\nu)}$ will be non-zero restricted to the mesh boundary; the effect of the corresponding control points on the boundary will be analogous to the effect of the control points of a univariate $C^1$ quadratic spline curve.

### 3.6. Properties of almost-$C^1$ splines

We collect the properties of $\boldsymbol{x}$ and $\mathscr{B}$ in the following results. In particular, the definitions of spline functions outlined in Section 3.5 immediately imply properties that are useful in numerical simulations.

**Lemma 3.6.** *For $\gamma \in \mathscr{T}_0^E$, consider the set of local control points $\{\boldsymbol{x}_\phi \;:\; \gamma \in \phi \in \mathscr{I}^*\}$. Let all control points in this set be coplanar, with the common plane defined by the normal $\boldsymbol{n}_\gamma$. Then, in Equation* (13)*, $\mathbf{c}^i_{jk} = \hat{\mathbf{c}}_{jk}[\boldsymbol{x}^*;\sigma_i]$.*

**Corollary 3.7.** *For all $\gamma \in \mathscr{T}_0^E$, let the set of local control points $\{\boldsymbol{x}_\phi \;:\; \gamma \in \phi \in \mathscr{I}^*\}$ be coplanar with the common plane defined by the normal $\boldsymbol{n}_\gamma$. Then, $\boldsymbol{x} = \boldsymbol{x}^*$.*

**Proposition 3.8.**

*(a) The total number of dofs satisfies*

$$\mathfrak{n} := |\mathscr{I}| \leq |\mathscr{T}_2| + |\mathscr{T}_1^B| + |\mathscr{T}_0^C| + 3|\mathscr{T}_0^E| \; . \tag{15}$$

*If no face and no boundary edge of the mesh contains only extraordinary vertices, then the above relation becomes an equality.*

*(b) Non-negativity: On any $\sigma \in \mathscr{T}_2$ and any $\phi \in \mathscr{I}$, $B_\phi\big|_\sigma \geq 0$.*

*(c) Partition of unity: On any $\sigma \in \mathscr{T}_2$, $\sum_{\phi \in \mathscr{I}} B_\phi\big|_\sigma \equiv 1$.*

*(d) Local support: If $B_\phi$ is such that $\phi \in \mathscr{T}_k$, $k = 0,1,2$, then $B_\phi\big|_\sigma = 0$ for any $\sigma \in \mathscr{T}_2$ such that $\sigma \cap \phi = \emptyset$. Similarly, if $\gamma \in \mathscr{T}_0^E$ with neighbouring faces $\{\sigma_1,\ldots,\sigma_\mu\}$, then $B_{(\gamma,k)}|_\sigma = 0$ for any $\sigma \in \mathscr{T}_2 \setminus \{\sigma_1,\ldots,\sigma_\mu\}$.*

*(e) Boundary Kronecker–Delta: All $B_\phi$ with $\phi \in \mathscr{T}_2$ are identically zero on the boundary of $\mathscr{T}$.*

*(f) Linear independence: $\{B_\phi : \phi \in \mathscr{I}\}$ form a basis for $\mathscr{B}$.*

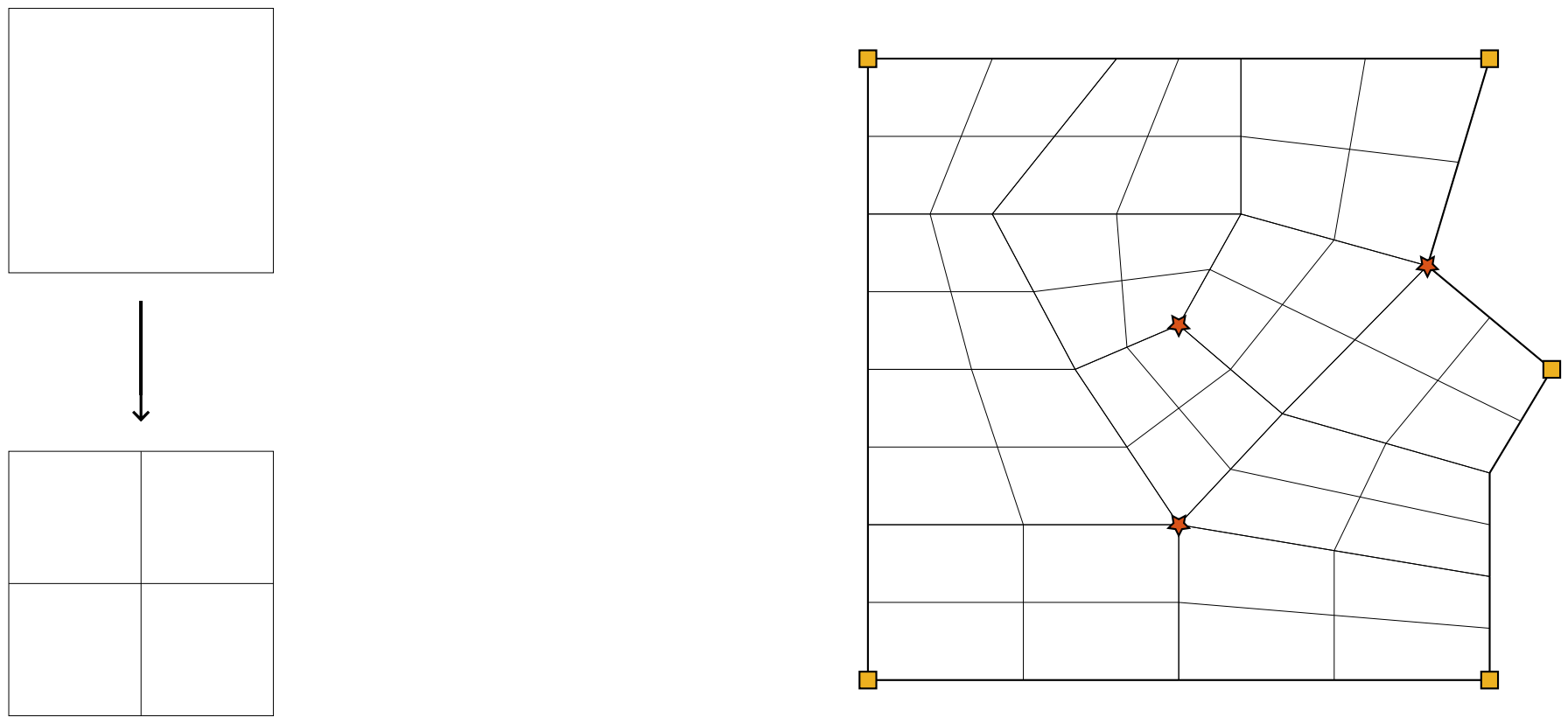

(a) refining a quadrilateral face

(b) refining the mesh from Figure 2(a)

Figure 7: The above figure shows how the mesh faces are split when performing refinement. Each quadrilateral face is split into 4 new faces as shown in (a). Figure (b) shows how the mesh topology changes after refinement for the mesh previously shown in Figure 2(a); it is also shown how the corner/extraordinary vertices retain their labels after refinement. See Section 4 for details.

In contrast to $\mathcal{B}^*$, the almost-$C^1$ splines spanning $\mathcal{B}$ are $C^1$ at all vertices, that is, any surface $\boldsymbol{x} \in \mathcal{B}^d$ possesses a well-defined tangent plane in every vertex, if the parameterization is regular. The tangent plane is orthogonal to the prescribed normal vector $\boldsymbol{n}_\gamma$ for each extraordinary vertex $\gamma \in \mathcal{T}_0^E$. Furthermore, for any spline surface $\boldsymbol{x} \in \mathcal{B}^d$ or $\boldsymbol{x}^* \in (\mathcal{B}^*)^d$ the boundary of the domain can be interpreted as a collection of quadratic B-spline curve segments with uniform (open) knot vectors. Similar to tensor-product B-splines, rational representations derived from $\mathcal{B}$ yield boundary curves that are quadratic NURBS. If the control points also satisfy $\boldsymbol{x}_\gamma = \boldsymbol{x}^*_\gamma$ for all $\gamma \in \mathcal{T}_0^C$, $\boldsymbol{x}_\tau = \boldsymbol{x}^*_\tau$ for all $\tau \in \mathcal{T}_1^B$, then the control points corresponding to boundary extraordinary vertices $\gamma \in \mathcal{T}_0^E \cap \mathcal{T}_0^B$ can be chosen such that the boundary curves of $\boldsymbol{x}$ and $\boldsymbol{x}^*$ coincide.

In Appendix C we present an alternative construction of almost-$C^1$ splines, which does not rely on subdividing the elements near extraordinary vertices. This results in a space that possesses no dofs corresponding to extraordinary faces, while three dofs per extraordinary vertex remain. The functions are biquadratic polynomials on all faces.

## 4. Refining almost-$C^1$ splines

Refinement of the mesh can help improve the resolving power of splines for the purpose of, for instance, obtaining a better approximation to the solution of a PDE. Assume that we are given a mesh $\mathcal{T}$, the associated almost-$C^1$ spline space $\mathcal{B}$, and control points $\boldsymbol{x}_\phi$, $\phi \in \mathcal{I}$, that define a spline geometry $\boldsymbol{x}$. In this section, we outline precisely how they can be refined. We start by describing the refinement of $\mathcal{T}$ in Section 4.1; next, in Section 4.2 we explain the motivation behind the refinement scheme for $\boldsymbol{x}$ and $\mathcal{B}$, the latter is described in Sections 4.3–4.5. Finally, we discuss some properties of the refinement scheme in Section 4.6. Our refinement scheme is closely related to the one presented recently in [64].

### *4.1. Refining $\mathcal{T}$*

To refine $\mathcal{T}$, which only contains topological information, we only need to specify how the connectivity and quadrilateral-composition of $\mathcal{T}$ are to be updated. In this document, we adopt a simple global refinement approach whereby all quadrilaterals are split into $2 \times 2$ quadrilaterals; see Figure 7. We will denote refined quantities with a “hat” — for instance, the refined mesh will be denoted as $\widehat{\mathcal{T}}$. During this process, no hanging nodes are introduced, i.e., $|\mathcal{T}_1| + |\mathcal{T}_2|$ new mesh vertices are added. We assume that the old mesh vertices retain their labels as per Section 2. Moreover, since the old vertices retain their labels, the corner and extraordinary vertices of $\widehat{\mathcal{T}}$ are respectively chosen to be identical to the corner and extraordinary vertices of $\mathcal{T}$, i.e., $\widehat{\mathcal{T}}_0^C = \mathcal{T}_0^C$ and $\widehat{\mathcal{T}}_0^E = \mathcal{T}_0^E$. Finally, the refined dof index sets $\widehat{\mathcal{I}}^*$ and $\widehat{\mathcal{I}}$ are defined for $\widehat{\mathcal{T}}$ exactly as outlined in Section 3. We directly obtain

$$\widehat{\mathcal{I}} = \widehat{\mathcal{T}}_2 \bigcup \widehat{\mathcal{T}_1^B} \bigcup \widehat{\mathcal{T}_0^C} \bigcup (\widehat{\mathcal{T}_0^E} \times \{1,2,3\}).$$

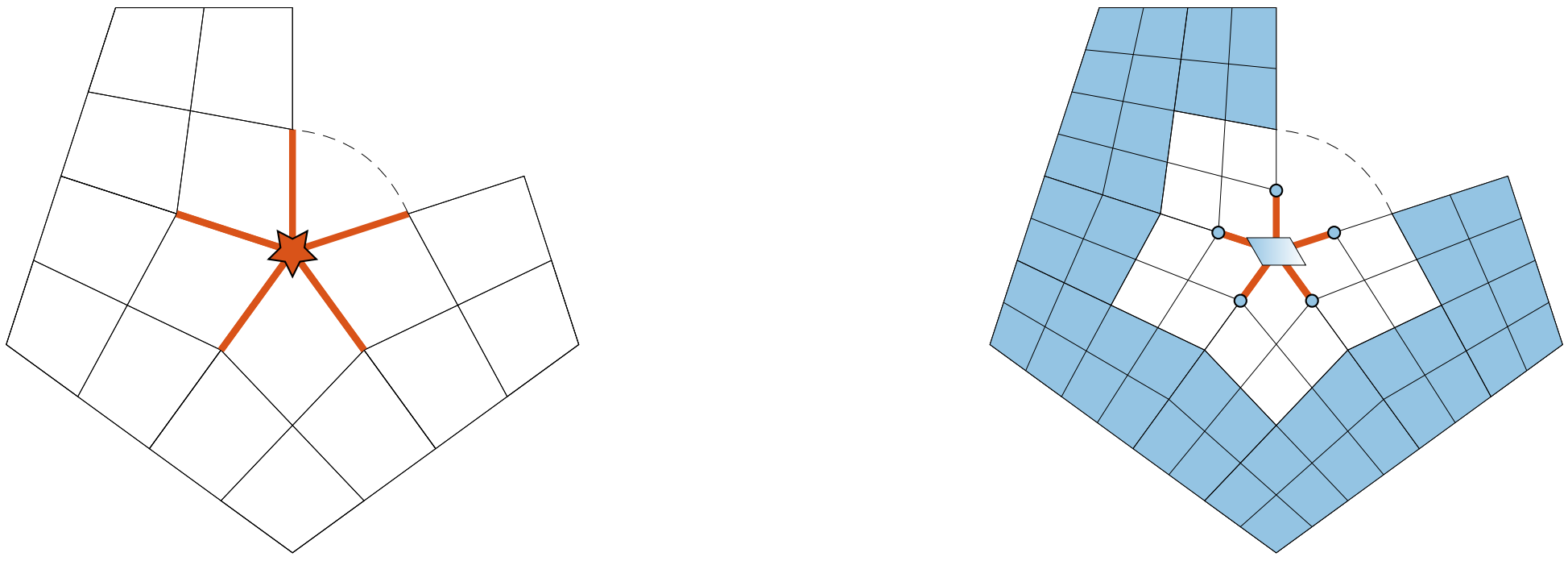

(a) coarse neighborhood (b) refined neighborhood

Figure 8: The above schematic depicts the motivation for our non-nested refinement scheme (c.f. Section 4.2 and Proposition 4.3). In (a) is a 2-neighbourhood of an extraordinary vertex; it is assumed that except the central vertex, all other vertices in this neighbourhood are regular. In (b) is the refinement of the coarse neighbourhood from (a). In both (a) and (b), splines (in $\mathscr{B}$ and $\widehat{\mathscr{B}}$, respectively) are not $C^1$ smooth across the *interiors* of the bold red edges. In (b), the flat shaded elements correspond to those where the refined spline geometry is identical to the coarse geometry; the filled disks correspond to points where the refined spline geometry interpolates the coarse geometry; and the central shaded parallelogram indicates that the coarse and refined spline geometries have identical normals at the extraordinary vertex.

### *4.2. Motivation behind refinement of $\boldsymbol{x}$ and $\mathscr{B}$*

In the next three subsections, we describe the refinement scheme in three steps. First, in Section 4.3, we use B-spline knot insertion to refine the element-local restrictions $\boldsymbol{x}|_\sigma$ on all elements $\sigma$ of $\mathscr{T}$. Let the thus refined element-local control points for $\sigma$ be indexed as $\widehat{\boldsymbol{c}}^{\sigma}_{jk}$, $0 \le j,k \le 3$. Next, in Section 4.4, we obtain the refined spline control points $\widehat{\boldsymbol{x}}^*_\phi$, $\phi \in \mathscr{I}^*$, as linear combinations of the control points $\widehat{\boldsymbol{c}}_{jk,\sigma}$. These define a geometry $\widehat{\boldsymbol{x}}^* \in (\widehat{\mathscr{B}}^*)^d$ on $\widehat{\mathscr{T}}$. Finally, in Section 4.5, we use $\widehat{\boldsymbol{x}}^*$ and the extraordinary-vertex normals $\boldsymbol{n}_\gamma$, $\gamma \in \widehat{\mathscr{T}}_0^E$, to define both $\widehat{\mathscr{B}}$ and $\widehat{\boldsymbol{x}}$. It is worth mentioning that, as shown in Section 4.6, $\widehat{\boldsymbol{x}}$ will be identical to $\widehat{\boldsymbol{x}}^*$.

Thus, the refinement scheme as outlined above will only ensure that $\widehat{\boldsymbol{x}} = \widehat{\boldsymbol{x}}^* \approx \boldsymbol{x}$ and, in general, the spline spaces will be non-nested, i.e., $\mathscr{B} \not\subset \widehat{\mathscr{B}}^*$ and $\mathscr{B} \not\subset \widehat{\mathscr{B}}$. We opt for such non-nested refinements because nested refinements would necessarily require involved bookkeeping (e.g., introduction of additional dofs with a refinement-level-dependent structure [63, 67]) which we seek to avoid.

The non-nested refinement scheme implies that the map $\boldsymbol{x} \mapsto \widehat{\boldsymbol{x}}$ can be specified in different ways. For instance, if there is a 'true geometry' $\boldsymbol{X}$ (e.g., the geometry at the coarsest refinement level or an underlying smooth surface), then at any given refinement level $\widehat{\boldsymbol{x}}$ can be computed so as to minimize $\|\widehat{\boldsymbol{x}} - \boldsymbol{X}\|^2$ in a suitable norm. Alternatively, the refined geometry $\widehat{\boldsymbol{x}}$ can be computed so that certain desirable properties of $\boldsymbol{x}$ are preserved – this is the approach we adopt. We formulate a refinement scheme that only uses local information and achieves two objectives. Firstly, it ensures that $\boldsymbol{x}$ and $\widehat{\boldsymbol{x}}$ are equal in the structured parts of the mesh $\mathscr{T}$ – in these parts $\boldsymbol{x}$ coincides with a tensor-product biquadratic spline, and thus B-spline knot insertion is sufficient for achieving this objective. Secondly, certain values and derivatives of $\boldsymbol{x}$ and $\widehat{\boldsymbol{x}}$ are equal in the unstructured parts of $\mathscr{T}$. See Figure 8 for a more precise visual overview of the scheme.

The refinement scheme will be depicted graphically in the following sections. The old control points will not be shown in the figures and, following our earlier convention, the new control points will be shown as filled blue circles. To further declutter the figures, only the indices of the new control points will be annotated.

### *4.3. Refining element-local representations $\boldsymbol{x}|_\sigma$*

On face $\sigma$ of $\mathscr{T}$, consider the restriction $\boldsymbol{x}|_\sigma$. If $\sigma \in \mathscr{T}_2 \backslash \mathscr{T}_2^E$, then $\boldsymbol{x}|_\sigma$ is a linear combination of the Bernstein–Bézier polynomials $b^0_{jk,\square}$ for some element-local control points $\boldsymbol{c}_{jk}[\boldsymbol{x};\sigma]$,

$$\boldsymbol{x}|_\sigma =: \sum_{j,k=0}^{2} \boldsymbol{c}_{jk}[\boldsymbol{x};\sigma] b^0_{jk,\square} \,. \tag{16}$$

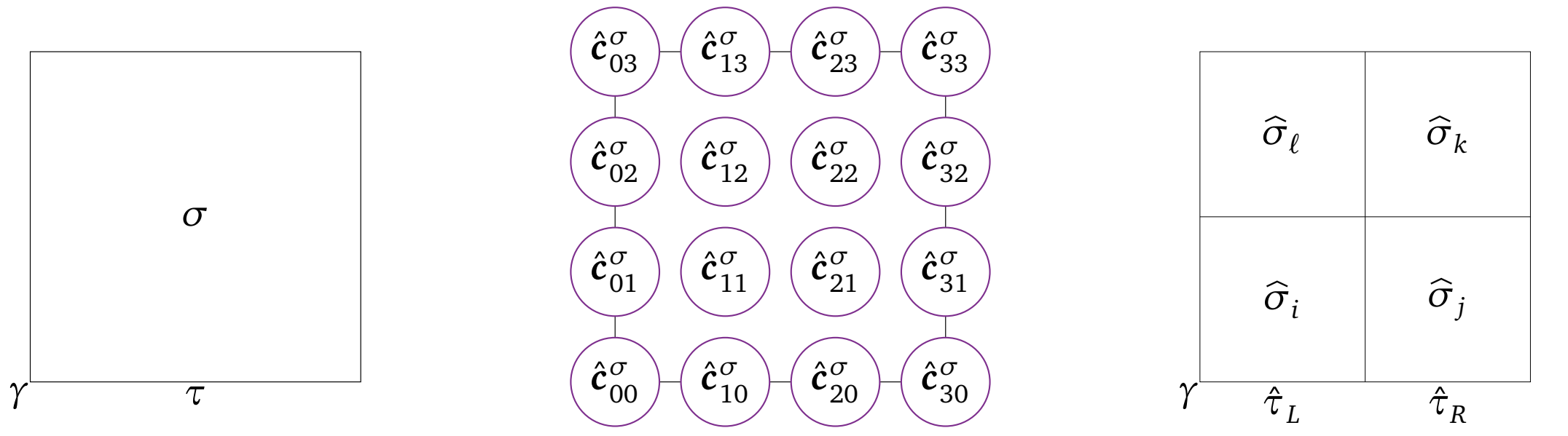

(a) coarse mesh face $\sigma$ (b) refined element-local control points (c) mesh faces after refining $\sigma$

Figure 9: The above figures show (a) a face $\sigma$ that belongs to the coarse mesh $\mathscr{T}$, (b) the refined element-local control points corresponding to the coarse spline geometry $\boldsymbol{x}|_\sigma$ (see Section 4.3), and (c) the faces of $\widehat{\mathscr{T}}$ obtained by refinements of $\sigma$. In particular, these figures are used as reference to explain the refinement scheme in Section 4.4.

Similarly, if $\sigma \in \mathscr{T}_2^E$, then $\boldsymbol{x}|_\sigma$ is a linear combination of the $C^1$ bi-quadratic B-splines $b^1_{jk,\square}$ for some element-local control points $\boldsymbol{c}_{jk}[\boldsymbol{x};\sigma]$,

$$\boldsymbol{x}|_\sigma =: \sum_{j,k=0}^{3} \boldsymbol{c}_{jk}[\boldsymbol{x};\sigma] b^1_{jk,\square} \,. \tag{17}$$

Then, for any $\sigma \in \mathscr{T}_2$ and with the matrix $\boldsymbol{c}[\boldsymbol{x};\sigma] := \big(\boldsymbol{c}_{jk}[\boldsymbol{x};\sigma]\big)_{jk}$, we define the element-local refined control points $\hat{\boldsymbol{c}}_{jk}[\boldsymbol{x};\sigma]$, $0 \le j,k \le 3$, as

$$\hat{\boldsymbol{c}}[\boldsymbol{x};\sigma] = \begin{cases} K\boldsymbol{c}[\boldsymbol{x};\sigma]K^T \,, & \sigma \in \mathscr{T}_2 \backslash \mathscr{T}_2^E \,, \\ \boldsymbol{c}[\boldsymbol{x};\sigma] \,, & \sigma \in \mathscr{T}_2^E \,, \end{cases} \tag{18}$$

where $K$ is the univariate B-spline knot-insertion matrix as in (11). Note that for all faces $\sigma$ in $\mathscr{T}_2$, we have

$$\boldsymbol{x}|_\sigma = \sum_{j,k=0}^{3} \hat{\boldsymbol{c}}_{jk}[\boldsymbol{x};\sigma] b^1_{jk,\square} \,. \tag{19}$$

For later reference, we will use the schematic shown in Figure 9 where we use the shorthand $\hat{\boldsymbol{c}}^\sigma_{jk} := \hat{\boldsymbol{c}}_{jk}[\boldsymbol{x};\sigma]$ for convenience.

### 4.4. *Defining $\widehat{\boldsymbol{x}}^*$ and $\widehat{\mathscr{B}}^*$*

The spline space $\widehat{\mathscr{B}}^*$ and the associated dof index set $\widehat{\mathscr{I}}^*$ are defined on $\widehat{\mathscr{T}}$ following the approach in Section 3.3. Then, we define a geometry $\widehat{\boldsymbol{x}}^* \in (\widehat{\mathscr{B}}^*)^d$ by computing the associated control points $\widehat{\boldsymbol{x}}^*_\phi$, $\phi \in \widehat{\mathscr{I}}^*$, as linear combinations of the element-local control points $\hat{\boldsymbol{c}}[\boldsymbol{x};\sigma]$, $\sigma \in \mathscr{T}_2$, from Section 4.3. We split this computation in four parts: corner vertex control points, boundary-edge control points, face control points in locally structured regions, and face control points in locally unstructured regions.

**Corner vertex control points**

As mentioned earlier, the number of corner vertices remains fixed during refinement (see Figure 7 for an example) and these vertices retain their labels. With reference to Figure 9(a) and (c), let $\gamma \in \widehat{\mathscr{T}}_0^C = \mathscr{T}_0^C$ and let $\widehat{\boldsymbol{x}}^*_\gamma$ be the refined spline control point associated to it. Then, we set

$$\widehat{\boldsymbol{x}}^*_\gamma := \hat{\boldsymbol{c}}_{00}[\boldsymbol{x};\sigma] \,. \tag{20}$$

**Boundary-edge control points**

With reference to Figure 9(a), let $\tau \in \mathscr{T}_1^B$ be refined into two new boundary edges $\widehat{\tau}_L, \widehat{\tau}_R \in \widehat{\mathscr{T}}_1^B$, as shown in Figure 9(c). Then, the control points corresponding to $\widehat{\tau}_L$ and $\widehat{\tau}_R$, respectively denoted as $\widehat{\boldsymbol{x}}^*_L$ and $\widehat{\boldsymbol{x}}^*_R$ here, are defined as,

$$\widehat{\boldsymbol{x}}^*_L := \hat{\boldsymbol{c}}_{10}[\boldsymbol{x};\sigma] \,, \quad \widehat{\boldsymbol{x}}^*_R := \hat{\boldsymbol{c}}_{20}[\boldsymbol{x};\sigma] \,.$$

**Face control points: locally structured regions**

A locally structured part of $\mathscr{T}$ is composed of boundary edges and faces that do not contain any extraordinary vertices. With reference to Figure 9(a), let $\sigma \in \mathscr{T}_2$ be a face containing no extraordinary vertices. Let $\sigma$ be refined into four new faces $\widehat{\sigma}_i, \widehat{\sigma}_j, \widehat{\sigma}_k$ and $\widehat{\sigma}_\ell$, as shown in Figure 9(c). Then, the corresponding control points, respectively denoted as $\widehat{\boldsymbol{x}}^*_i, \widehat{\boldsymbol{x}}^*_j, \widehat{\boldsymbol{x}}^*_k$ and $\widehat{\boldsymbol{x}}^*_\ell$, are defined as,

$$\left(\widehat{\boldsymbol{x}}^*_i, \widehat{\boldsymbol{x}}^*_j, \widehat{\boldsymbol{x}}^*_k, \widehat{\boldsymbol{x}}^*_\ell\right) := \left(\hat{\boldsymbol{c}}_{11}[\boldsymbol{x};\sigma],\ \hat{\boldsymbol{c}}_{21}[\boldsymbol{x};\sigma],\ \hat{\boldsymbol{c}}_{22}[\boldsymbol{x};\sigma],\ \hat{\boldsymbol{c}}_{12}[\boldsymbol{x};\sigma]\right).$$

**Face control points: locally unstructured regions**

An unstructured part of $\mathscr{T}$ is a face that contains one or more extraordinary vertices. Since we opt for non-nested refinements, it is not possible to exactly preserve the geometry in the unstructured parts of $\mathscr{T}$ during refinement. Nevertheless, the following approach ensures that $\widehat{\boldsymbol{x}}^*$ interpolates the midpoints of the coarse spoke edges when all faces of $\mathscr{T}$ contain at most one extraordinary vertex.

First, with reference to Figure 9(a), let $\sigma \in \mathscr{T}_2$ be a face containing either more than one extraordinary vertices, or containing an extraordinary vertex that is in the 1-ring of another extraordinary vertex. This $\sigma$ is being refined into four new faces $\widehat{\sigma}_i, \widehat{\sigma}_j, \widehat{\sigma}_k$ and $\widehat{\sigma}_\ell$, as shown in Figure 9(c). Then, similarly to the locally structured case, the corresponding control points, respectively denoted as $\widehat{\boldsymbol{x}}^*_i, \widehat{\boldsymbol{x}}^*_j, \widehat{\boldsymbol{x}}^*_k$ and $\widehat{\boldsymbol{x}}^*_\ell$, are defined as

$$\left(\widehat{\boldsymbol{x}}^*_i, \widehat{\boldsymbol{x}}^*_j, \widehat{\boldsymbol{x}}^*_k, \widehat{\boldsymbol{x}}^*_\ell\right) := \left(\hat{\boldsymbol{c}}_{11}[\boldsymbol{x};\sigma],\ \hat{\boldsymbol{c}}_{21}[\boldsymbol{x};\sigma],\ \hat{\boldsymbol{c}}_{22}[\boldsymbol{x};\sigma],\ \hat{\boldsymbol{c}}_{12}[\boldsymbol{x};\sigma]\right).$$

Now we tackle the final remaining case: refinement of faces around an extraordinary vertex that is not contained in the 1-ring of any other extraordinary vertex. We first consider the case where this extraordinary vertex is an interior vertex and later when it is a boundary vertex.

Let $\gamma \in \mathscr{T}_0^E \cap \mathring{\mathscr{T}}_0$ be an interior extraordinary vertex of valence $\mu$ and let the labelling of the coarse faces and their refined element-local control points in the neighborhood of $\gamma$ be as in Figure 10(a). Then, the corresponding refined control points $\widehat{\boldsymbol{x}}^*_{r,i}, \widehat{\boldsymbol{x}}^*_{r,j}, \widehat{\boldsymbol{x}}^*_{r,k}$ and $\widehat{\boldsymbol{x}}^*_{r,\ell}$, $r = 1, \dots, \mu$, are computed as below, where we employ the shorthand $\widehat{\boldsymbol{c}}^r_{jk} := \hat{\boldsymbol{c}}_{jk}[\boldsymbol{x};\sigma_r]$,

$$\begin{aligned}
&\begin{bmatrix} \widehat{\boldsymbol{x}}^*_{r,j} \\ \widehat{\boldsymbol{x}}^*_{r,k} \\ \widehat{\boldsymbol{x}}^*_{r,\ell} \end{bmatrix} = \begin{bmatrix} \widehat{\boldsymbol{c}}^r_{21} \\ \widehat{\boldsymbol{c}}^r_{22} \\ \widehat{\boldsymbol{c}}^r_{12} \end{bmatrix}, \quad r = 1, \dots, \mu, \\
&\begin{bmatrix} \widehat{\boldsymbol{x}}^*_{1,i} \\ \widehat{\boldsymbol{x}}^*_{2,i} \\ \vdots \\ \widehat{\boldsymbol{x}}^*_{\mu,i} \end{bmatrix} = \mathring{\mathrm{S}}_\mu \begin{bmatrix} \widehat{\boldsymbol{c}}^1_{01} \\ \widehat{\boldsymbol{c}}^2_{01} \\ \vdots \\ \widehat{\boldsymbol{c}}^\mu_{01} \end{bmatrix} + \mathring{\mathrm{Q}}_\mu \begin{bmatrix} \widehat{\boldsymbol{c}}^1_{11} \\ \widehat{\boldsymbol{c}}^2_{11} \\ \vdots \\ \widehat{\boldsymbol{c}}^\mu_{11} \end{bmatrix},
\end{aligned} \tag{21}$$

where $\mathring{\mathrm{S}}_\mu$ and $\mathring{\mathrm{Q}}_\mu$ are circulant matrices. If $\mu$ is odd, they are defined to be

$$\mathring{\mathrm{S}}_\mu = \mathrm{circulant}\big(1,\ -1,\ 1,\ -1,\ \cdots,\ -1,\ 1\big), \quad \mathring{\mathrm{Q}}_\mu = \mathbf{0}, \tag{22}$$

and if $\mu$ is even, they are defined to be

$$\begin{aligned}
\mathring{\mathrm{S}}_\mu &= \mathrm{circulant}\left(2 - \frac{2}{\mu},\ -\left(2 - \frac{4}{\mu}\right),\ 2 - \frac{6}{\mu},\ -\left(2 - \frac{8}{\mu}\right),\ \cdots,\ 0\right), \\
\mathring{\mathrm{Q}}_\mu &= \frac{1}{\mu}\mathrm{circulant}(1,\ -1,\ 1,\ -1,\ \cdots,\ -1).
\end{aligned} \tag{23}$$

Next, let $\gamma \in \mathscr{T}_0^E \cap \mathscr{T}_0^B$ be a boundary extraordinary vertex of valence $\mu$, and let the labelling of the coarse faces and their refined element-local control points in the neighbourhood of $\gamma$ be as in Figure 10(b). Then, the corresponding refined control points $\widehat{\boldsymbol{x}}^*_{r,i}, \widehat{\boldsymbol{x}}^*_{r,j}, \widehat{\boldsymbol{x}}^*_{r,k}$ and $\widehat{\boldsymbol{x}}^*_{r,\ell}$, $r = 1, \dots, \mu$, are computed as below, where we again employ the shorthand $\widehat{\boldsymbol{c}}^r_{jk} := \hat{\boldsymbol{c}}_{jk}[\boldsymbol{x};\sigma_r]$,

$$\begin{bmatrix}\widehat{\boldsymbol{x}}^*_{r,j}\\ \widehat{\boldsymbol{x}}^*_{r,k}\\ \widehat{\boldsymbol{x}}^*_{r,\ell}\end{bmatrix} = \begin{bmatrix}\widehat{\boldsymbol{c}}^r_{21}\\ \widehat{\boldsymbol{c}}^r_{22}\\ \widehat{\boldsymbol{c}}^r_{12}\end{bmatrix}, \quad r=1,\dots,\mu\,,$$

$$\begin{bmatrix}\widehat{\boldsymbol{x}}^*_{1,i}\\ \widehat{\boldsymbol{x}}^*_{2,i}\\ \vdots\\ \widehat{\boldsymbol{x}}^*_{\mu,i}\end{bmatrix} = \mathrm{S}^\partial_\mu \begin{bmatrix}\widehat{\boldsymbol{c}}^1_{01}\\ \widehat{\boldsymbol{c}}^2_{01}\\ \vdots\\ \widehat{\boldsymbol{c}}^{\mu-1}_{01}\end{bmatrix} + \mathrm{Q}^\partial_\mu \begin{bmatrix}\widehat{\boldsymbol{c}}^1_{11}\\ \widehat{\boldsymbol{c}}^2_{11}\\ \vdots\\ \widehat{\boldsymbol{c}}^{\mu}_{11}\end{bmatrix}, \tag{24}$$

where $\mathrm{S}^\partial_\mu$ is defined to be

$$\mathrm{S}^\partial_\mu = \frac{\mathrm{R}_\mu + \mathrm{J}_\mu \mathrm{R}_\mu \mathrm{J}_{\mu-1}}{2}\,, \tag{25}$$

and $\mathrm{Q}^\partial_\mu$, $\mathrm{R}_\mu$ and $\mathrm{J}_k$ are defined to be the following matrices of sizes $\mu\times\mu$, $\mu\times(\mu-1)$ and $k\times k$, $k\geq 1$, respectively,

$$\mathrm{Q}^\partial_\mu = \frac{1}{\mu}\left((-1)^{j+k}\right)_{jk}\,, \quad \mathrm{R}_\mu = \begin{bmatrix} 4-\frac{4}{\mu} & -\left(4-\frac{8}{\mu}\right) & \left(4-\frac{12}{\mu}\right) & \cdots & (-1)^{\mu}\frac{4}{\mu}\\ & \left(4-\frac{8}{\mu}\right) & -\left(4-\frac{12}{\mu}\right) & \cdots & (-1)^{\mu+1}\frac{4}{\mu}\\ & & \left(4-\frac{12}{\mu}\right) & \cdots & (-1)^{\mu+2}\frac{4}{\mu}\\ & & & \ddots & \vdots\\ & & & & (-1)^{2\mu-2}\frac{4}{\mu}\\ & & & & 0\end{bmatrix},$$

$$\mathrm{J}_k = \begin{bmatrix} & & & 1\\ & & \iddots & \\ & 1 & & \\ 1 & & & \end{bmatrix}.$$

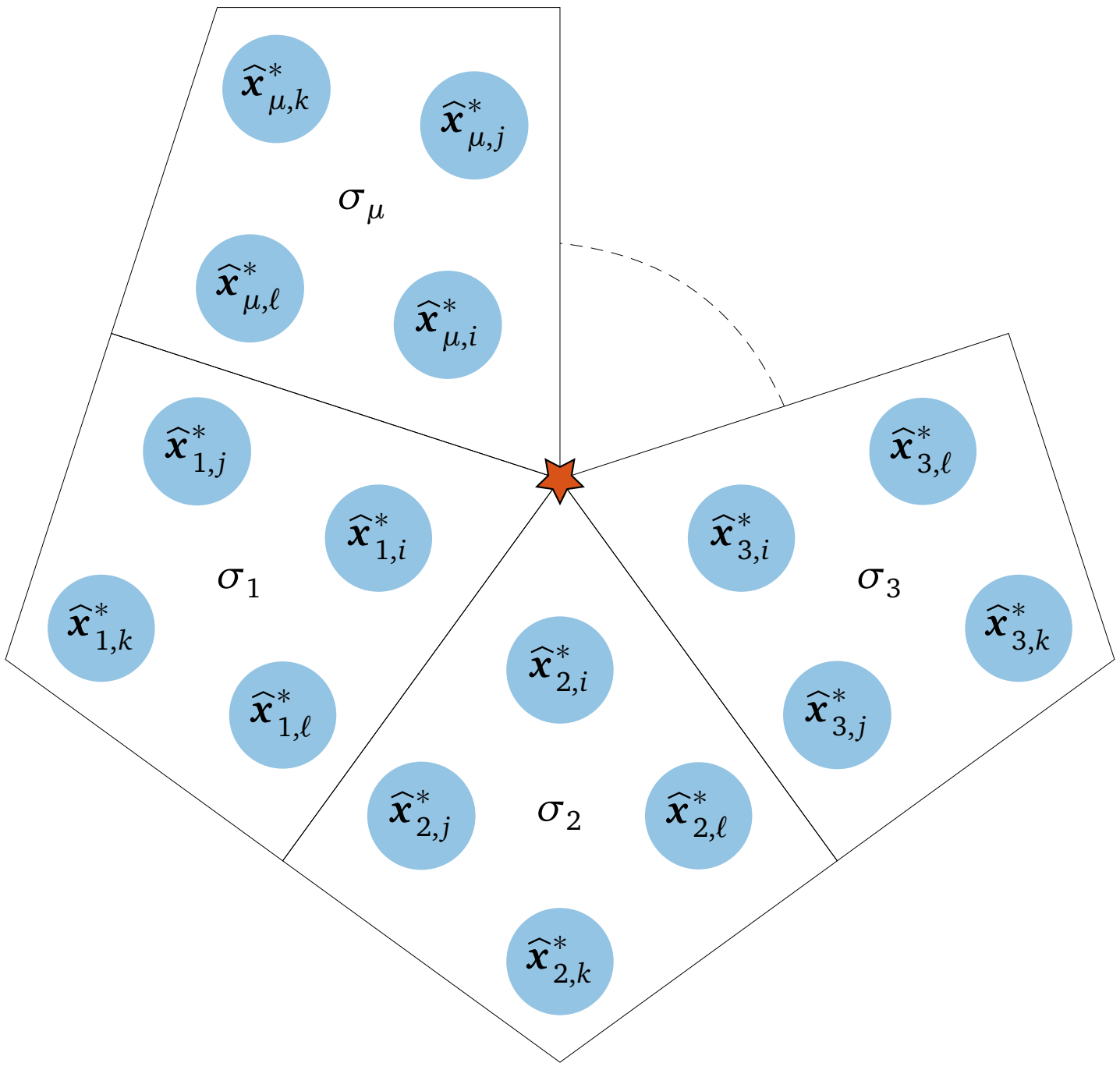

(a) interior extraordinary vertex

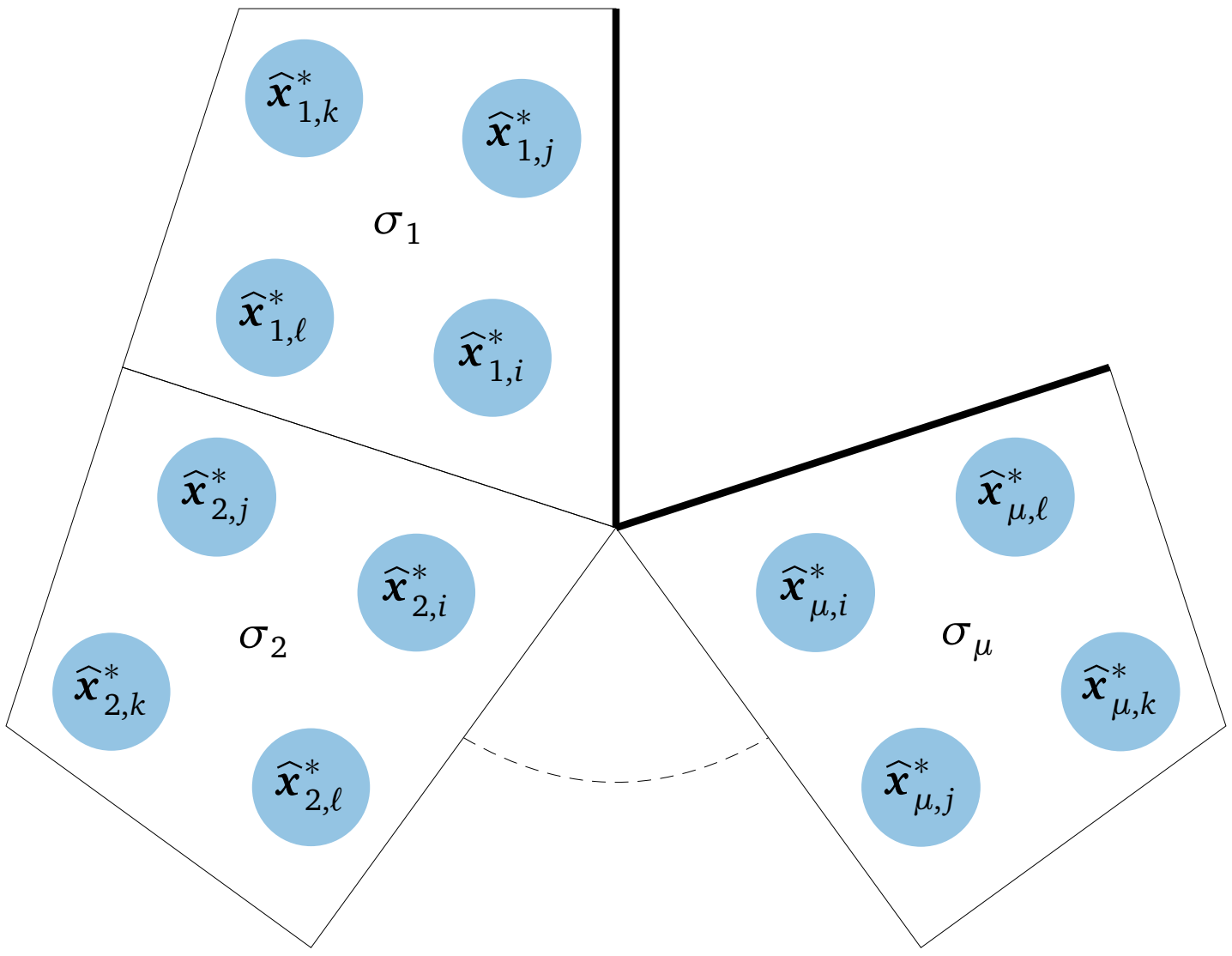

(b) boundary extraordinary vertex

Figure 10: The above shows the refinement of face control points in the locally unstructured regions of the mesh. Let $\sigma_1, \ldots, \sigma_\mu \in \mathscr{T}_2$ be quadrilaterals that share an extraordinary vertex and, moreover, this is the only extraordinary vertex that each face contains; figure (a) shows the case for an interior extraordinary vertex and figure (b) for a boundary extraordinary vertex. Then, the refined control points $\widehat{\boldsymbol{x}}^*_{m,i}$, $\widehat{\boldsymbol{x}}^*_{m,j}$, $\widehat{\boldsymbol{x}}^*_{m,k}$ and $\widehat{\boldsymbol{x}}^*_{m,\ell}$ are obtained using Equations (21), (22) and (23) for figure (a), and Equations (24) and (25) for figure (b). Note that the labelling of face control points with the subscripts $i, j, k, \ell$ corresponds to the refined face labelling shown in Figure 9(c), assuming that the extraordinary vertex in the above figures coincides with $\gamma$ in Figure 9(a).

*4.5. Defining the refined almost-$C^1$ spline geometry $\widehat{\boldsymbol{x}}$ and $\widehat{\mathscr{B}}$*

Finally, starting from $\widehat{\boldsymbol{x}}^*$ and the extraordinary-vertex normals $\boldsymbol{n}_\gamma$, $\gamma\in\widehat{\mathscr{T}}^E_0 = \mathscr{T}^E_0$, we define both $\widehat{\mathscr{B}}$ and $\widehat{\boldsymbol{x}}$. This is done exactly as outlined in Section 3.5.

*4.6. Properties of the refinement scheme*

The following results outline the properties of the refinement scheme explained in the previous sections. These properties were previously visually depicted in Figure 8.

**Lemma 4.1.** *The control points in the set $\left\{\widehat{\boldsymbol{x}}^*_{r,i}\ :\ r=1,\dots,\mu\right\}$ defined by Equations* (21) *and* (24) *are coplanar, with the common plane defined by the normal vector $\boldsymbol{n}_\gamma$.*

*Proof.* By definition of $\boldsymbol{x}$ and $\hat{\boldsymbol{c}}[\boldsymbol{x};\sigma]$, $\sigma\in\mathscr{T}_2$, the element-local control points

$$\left\{\hat{\boldsymbol{c}}_{jk}[\boldsymbol{x};\sigma_r]\ :\ (j,k)\in\{(0,0),(1,0),(0,1),(1,1)\}\,,\quad r=1,\dots,\mu\right\},$$

are coplanar. More precisely, the plane in which they lie is orthogonal to the normal vector $\boldsymbol{n}_\gamma$ specified for the extraordinary vertex in Figure 10.

Then, first considering the interior vertex case, and observing that each row of $\mathring{\mathrm{S}}_\mu + \mathring{\mathrm{Q}}_\mu$ sums to 1, we have

$$\begin{bmatrix}\widehat{\boldsymbol{x}}^*_{1,i}-\widehat{\boldsymbol{c}}^1_{00}\\ \widehat{\boldsymbol{x}}^*_{2,i}-\widehat{\boldsymbol{c}}^1_{00}\\ \vdots\\ \widehat{\boldsymbol{x}}^*_{\mu,i}-\widehat{\boldsymbol{c}}^1_{00}\end{bmatrix} = \mathring{\mathrm{S}}_\mu \begin{bmatrix}\widehat{\boldsymbol{c}}^1_{01}-\widehat{\boldsymbol{c}}^1_{00}\\ \widehat{\boldsymbol{c}}^2_{01}-\widehat{\boldsymbol{c}}^1_{00}\\ \vdots\\ \widehat{\boldsymbol{c}}^{\mu}_{01}-\widehat{\boldsymbol{c}}^1_{00}\end{bmatrix} + \mathring{\mathrm{Q}}_\mu \begin{bmatrix}\widehat{\boldsymbol{c}}^1_{11}-\widehat{\boldsymbol{c}}^1_{00}\\ \widehat{\boldsymbol{c}}^2_{11}-\widehat{\boldsymbol{c}}^1_{00}\\ \vdots\\ \widehat{\boldsymbol{c}}^{\mu}_{11}-\widehat{\boldsymbol{c}}^1_{00}\end{bmatrix}.$$

Thus the vectors on the left hand side and the vectors on the right hand side are coplanar and the claim follows. The boundary vertex case can be similarly shown since each row of $\mathrm{S}^\partial_\mu + \mathrm{Q}^\partial_\mu$ also sums to 1. ■

**Corollary 4.2.** *The geometry $\widehat{\boldsymbol{x}}$ is identical to $\widehat{\boldsymbol{x}}^*$, i.e., $\widehat{\boldsymbol{x}} = \widehat{\boldsymbol{x}}^* \in (\widehat{\mathscr{B}} \cap \widehat{\mathscr{B}^*})^d$.*

*Proof.* The claim follows from Corollary 3.7 and Lemma 4.1. ∎

**Proposition 4.3.** *Let $\widehat{\boldsymbol{x}}$ be obtained from $\boldsymbol{x}$ via the refinement process outlined in Sections 4.3–4.5. Then the following hold true.*

*(a) Translation and rotation invariance: Let* $\mathrm{A}$ *be the matrix mapping the coarse control points $\boldsymbol{x}_\phi$ to the refined control points $\widehat{\boldsymbol{x}}_\phi$. If* $\mathrm{T}$ *and* $\mathrm{R}$ *denote a translation and a rotation, then*

$$\mathrm{A} \circ \mathrm{T} = \mathrm{T} \circ \mathrm{A}\,, \quad \mathrm{A} \circ \mathrm{R} = \mathrm{R} \circ \mathrm{A}\,.$$

*In particular, the rows of* $\mathrm{A}$ *sum to* $1$.

*(b) Boundary preservation: Let $\tau$ be a boundary edge of $\mathscr{T}$ and $\boldsymbol{x}_\tau$ the local representation of $\boldsymbol{x}$ restricted to $\tau$. With reference to Figure 9, let the origin of local coordinates $[0,1]$ on $\tau$ be at the left end. Then,*

$$\widehat{\boldsymbol{x}}|_{\widehat{\tau}_L} = \boldsymbol{x}_\tau|_{[0,0.5]}\,, \quad \widehat{\boldsymbol{x}}|_{\widehat{\tau}_R} = \boldsymbol{x}_\tau|_{[0.5,1]}\,.$$

*(c) Structured quadrilateral preservation: Let $\sigma$ be a quadrilateral of $\mathscr{T}$ that contains no extraordinary vertices and no corner vertices of valence $> 1$, and $\boldsymbol{x}_\sigma$ the local polynomial representation of $\boldsymbol{x}$ restricted to $\sigma$. With reference to Figure 9, let the origin of local coordinates $[0,1]^2$ on $\sigma$ be at the bottom-left vertex $\gamma$. Then,*

$$\widehat{\boldsymbol{x}}|_{\widehat{\sigma}_i} = \boldsymbol{x}_\sigma|_{[0,0.5]^2}\,, \quad \widehat{\boldsymbol{x}}|_{\widehat{\sigma}_j} = \boldsymbol{x}_\sigma|_{[0.5,1]\times[0,0.5]}\,, \quad \widehat{\boldsymbol{x}}|_{\widehat{\sigma}_k} = \boldsymbol{x}_\sigma|_{[0.5,1]^2}\,, \quad \widehat{\boldsymbol{x}}|_{\widehat{\sigma}_l} = \boldsymbol{x}_\sigma|_{[0,0.5]\times[0.5,1]}\,.$$

*(d) Midpoint interpolation on interior spoke edges: Consider the settings shown in Figure 10(a) and (b) with the extraordinary vertex $\gamma$. Let $\tau$ be an interior spoke edge in the figures such that $\gamma$ corresponds to its first endpoint, and let it be refined into edges $\widehat{\tau}_a$ and $\widehat{\tau}_b$ such that $\gamma$ now corresponds to the first endpoint of $\widehat{\tau}_a$. Denote as $\boldsymbol{x}_\tau$ the local representation of $\boldsymbol{x}$ restricted to $\tau$. Then,*

$$\widehat{\boldsymbol{x}}|_{\widehat{\tau}_a}(1) = \boldsymbol{x}_\tau(0.5)\,.$$

*(e) Normal vector preservation: The refined geometry $\widehat{\boldsymbol{x}}$ has well-defined normal vectors at all vertices, and these are identical to the normal vectors for the coarse geometry $\boldsymbol{x}$.*

*(f) Dimension of the refined spline space: We have $\widehat{\mathscr{I}} = \widehat{\mathscr{I}^*} \cup (\widehat{\mathscr{T}}_0^E \times \{1,2,3\})$ and*

$$\widehat{\mathrm{n}} := |\widehat{\mathscr{I}}| = |\widehat{\mathscr{T}}_2| + |\widehat{\mathscr{T}}_1^B| + |\widehat{\mathscr{T}}_0^C| + 3|\widehat{\mathscr{T}}_0^E|\,. \tag{26}$$

*(g) Convergent refinement scheme: The limit surface exists and has a well-defined tangent plane at all points.*

*Proof.* It is easy to see that (a), (e) and (f) hold by construction. Following the same line of reasoning adopted in [64, Proposition 4.2], it can be seen that properties (b)–(d) hold for the geometry $\widehat{\boldsymbol{x}}^*$. Then, from Corollary 4.2, properties (b)–(d) also hold for $\widehat{\boldsymbol{x}}$. Finally, from [64, Appendix], the limit surface exists when the refinement scheme is applied to geometries in $(\mathscr{B}^*)^d$. Thus, from Corollary 4.2, the same is true for geometries in $\mathscr{B}^d$. While this limit surface may not have a well-defined tangent plane at extraordinary points for all geometries in $(\mathscr{B}^*)^d$, it will be smooth there for geometries in $\mathscr{B}^d$ by virtue of property (e). This implies property (g) and completes the proof. ∎

## 5. Numerical tests

In this section we test the approximation properties and conditioning of almost-$C^1$ splines for several second and fourth order PDE model problems of practical relevance. Convergence tests and condition number growths for the Poisson and Biharmonic problems are presented in Section 5.1, the Scordelis–Lo thin shell benchmark in Section 5.2, the surface Cahn–Hilliard problem in Section 5.3 and the surface Laplace–Beltrami eigenvalue problem in Section 5.4.

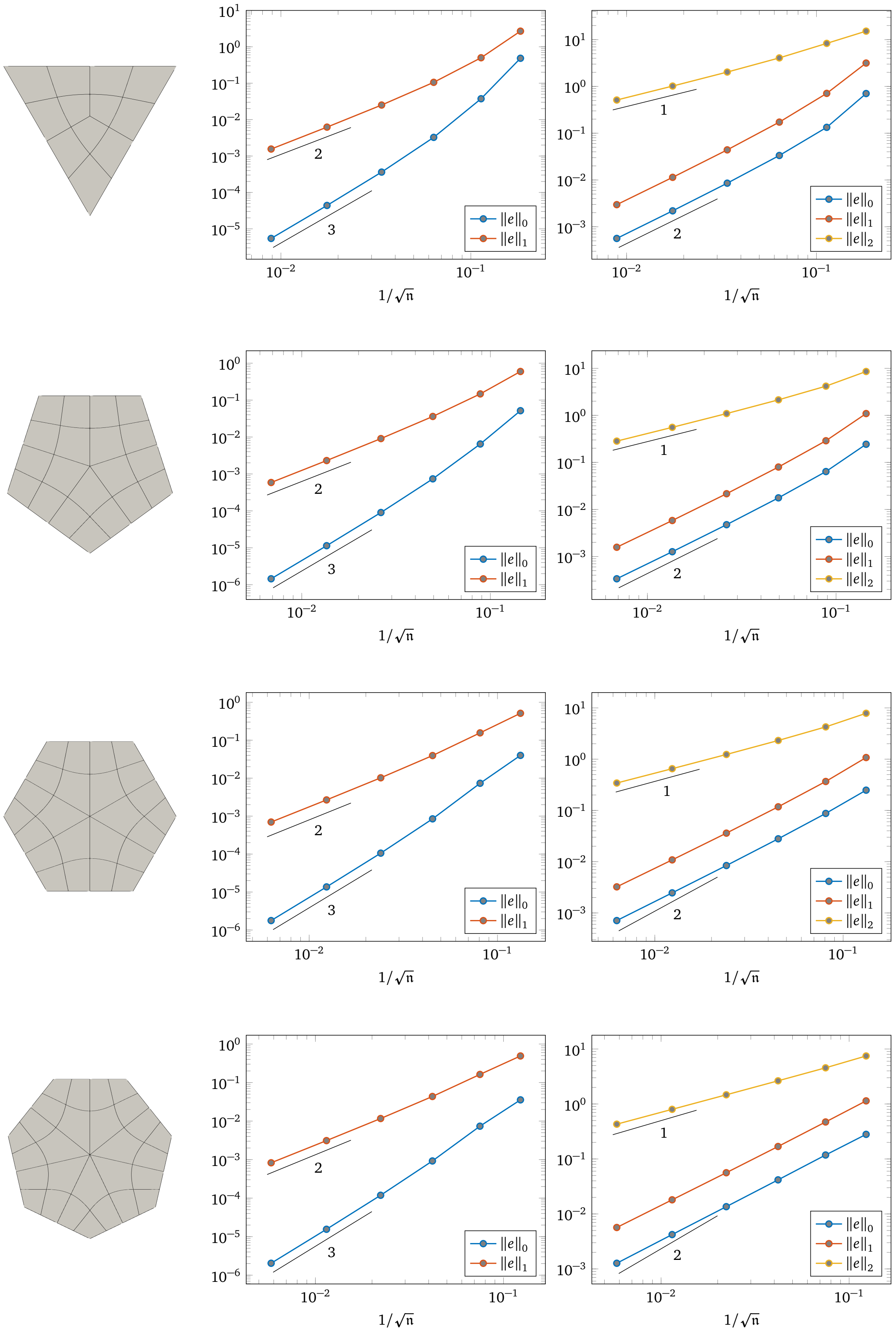

Figure 11: The above plots show the error convergence with mesh refinement when solving the Poisson (middle) and Biharmonic (right) problems in Equations (28) and (29), respectively, on the different spline geometries (left).

### 5.1. Poisson and Biharmonic problems

We start by investigating the convergence behaviour of almost-$C^1$ splines under mesh refinements. For this purpose, we solve $H^1$ and $H^2$ projection problems on meshes containing a single extraordinary vertex. The exact solution of the problem for all $\Omega$ is chosen to be

$$f_{\text{exact}}(x,y) = \sin\left(\pi x + \frac{\pi}{3}\right)\sin\left(\pi y + \frac{\pi}{5}\right) . \tag{27}$$

The following trial and test function spaces were used for the problem,

$$\mathscr{S}_0 := \left\{f \in \mathscr{B} : f = f_0 \text{ on } \Gamma\right\} , \qquad \mathscr{V}_0 := \left\{w \in \mathscr{B} : w = 0 \text{ on } \Gamma\right\} ,$$
$$\mathscr{S}_1 := \left\{f \in \mathscr{S}_0 : f_{,\boldsymbol{n}} = f_{1,\boldsymbol{n}} \text{ on } \Gamma\right\} , \qquad \mathscr{V}_1 := \left\{w \in \mathscr{V}_0 : w_{,\boldsymbol{n}} = 0 \text{ on } \Gamma\right\} ,$$

where $\Gamma = \partial\Omega$ is the boundary of $\Omega$, $\boldsymbol{n}$ is the unit normal to $\Gamma$, $a_{,\boldsymbol{n}} := \nabla a \cdot \boldsymbol{n}$, and $f_0 \in \mathscr{B}$ and $f_1 \in \mathscr{S}_0$ are defined such that,

$$\int_\Gamma w f_0 \, d\Gamma = \int_\Gamma w f_{\text{exact}} \, d\Gamma , \qquad \forall w \in \mathscr{B},\ w \neq 0 \text{ on } \Gamma ,$$
$$\int_\Gamma w_{,\boldsymbol{n}} f_{1,\boldsymbol{n}} \, d\Gamma = \int_\Gamma w_{,\boldsymbol{n}} f_{\text{exact},\boldsymbol{n}} \, d\Gamma , \qquad \forall w \in \mathscr{B},\ w_{,\boldsymbol{n}} \neq 0 \text{ on } \Gamma .$$

The weak form of the problems using the above trial and test spaces is

$$\mathrm{P}_1 : \text{Find } f \in \mathscr{S}_0 : \qquad \int_\Omega \nabla w \cdot \nabla f \, d\Omega = -\int_\Omega w \Delta f_{\text{exact}} \, d\Omega , \qquad \forall w \in \mathscr{V}_0 , \tag{28}$$

$$\mathrm{P}_2 : \text{Find } f \in \mathscr{S}_1 : \qquad \int_\Omega \Delta w \Delta f \, d\Omega = \int_\Omega w \Delta^2 f_{\text{exact}} \, d\Omega , \qquad \forall w \in \mathscr{V}_1 . \tag{29}$$

For different valent extraordinary points, the spline geometries $\Omega$ at the coarsest refinement level are shown on the left in Figure 11. The error convergence with mesh refinement is shown in the plots on the right, with the error norms plotted against the inverse of the square root of the number of dofs. We plot three different errors: $\|e\|_0$, $\|e\|_1$ and $\|e\|_2$, the $L^2(\Omega)$, $H^1(\Omega)$ and $H^2(\Omega)$ norms of the error, respectively. As is clear from Figure 11, the spline spaces demonstrate optimal convergence rates in the $L^2$, $H^1$ and $H^2$ norms. Only a slight deterioration is visible in the $L^2$ norm ($\approx 1.8$) for the Biharmonic problem at high refinement levels; this is similar to the results for bi-cubic splines in [63]. Finally, for both problems we show the condition numbers of the system matrices in Table 1, as can be seen, the conditioning is as expected from a well-conditioned basis without singularities.

*Remark* 5.1. Imposition of boundary conditions for values of unstructured splines is very simple here, given the boundary Kronecker–Delta property from Proposition 3.8. Boundary conditions for normal derivatives are also not problematic as they only involve the basis functions for the boundary dofs and the layer of dofs adjacent to boundary dofs (i.e., the adjacent face dofs and dofs for boundary extraordinary vertices).

### 5.2. Kirchhoff–Love shells: The Scordelis–Lo benchmark problem

We now solve a Kirchhoff–Love benchmark problem – the Scordelis–Lo test case. A curved cylinder with dimensions $(r,\ L,\ \theta_{\text{sector}}) = (25,\ 50,\ 2\pi/3)$ is loaded under gravity and has the following material parameters:

$$\text{Young's modulus,} \quad E = 4.32 \times 10^8 ,$$
$$\text{Poisson's ratio,} \quad \nu = 0.0 ,$$
$$\text{Thickness,} \quad t = 0.25 .$$

The shell formulation used is based on the Kirchhoff–Love thin shell theory in which transverse shear strains are zero. The end result is a rotation-free formulation requiring $C^1$ continuous trial functions. Necessarily, since transverse shear strains are suppressed, one would anticipate the theory would result in smaller deformations than for Reissner–Mindlin shell theory, which accounts for transverse shear strains. As observed in [33] and elsewhere, the Kirchhoff–Love theory leads to a converged maximum downward vertical displacement $w_{\text{ref}} = 0.3006$, while the Reissner–Mindlin theory yields 0.3024; in the following we use the former as the reference solution.

| **Problem** | $\mu$ | $\kappa_0$ | $\kappa_1$ | $\kappa_2$ | $\kappa_3$ | $\kappa_4$ | $\kappa_5$ | |
|---|---|---|---|---|---|---|---|---|
| $P_1$ | 3 | 5.10e+0 | 1.56e+1 | 5.51e+1 | 2.16e+2 | 8.78e+2 | 3.55e+3 | $\propto \mathfrak{n}$ |
| | 5 | 1.19e+1 | 2.30e+1 | 7.99e+1 | 3.04e+2 | 1.20e+3 | 4.85e+3 | |
| | 6 | 1.20e+1 | 2.25e+1 | 8.77e+1 | 3.50e+2 | 1.40e+3 | 5.69e+3 | |
| | 7 | 1.35e+1 | 2.89e+1 | 1.12e+2 | 4.51e+2 | 1.82e+3 | 7.44e+3 | |
| $P_2$ | 3 | 2.94e+1 | 6.23e+2 | 9.64e+3 | 1.55e+5 | 2.51e+6 | 4.09e+7 | $\propto \mathfrak{n}^2$ |
| | 5 | 3.63e+1 | 5.64e+2 | 8.46e+3 | 1.32e+5 | 2.10e+6 | 3.35e+7 | |
| | 6 | 3.39e+1 | 5.44e+2 | 8.27e+3 | 1.29e+5 | 2.07e+6 | 3.37e+7 | |
| | 7 | 4.32e+1 | 7.30e+2 | 1.13e+4 | 1.80e+5 | 2.95e+6 | 4.80e+7 | |

Table 1: The above table shows the condition numbers ($\kappa$) of the system matrices corresponding to problems $P_1$ and $P_2$ from Equation (28) and (29), respectively. The condition number at the $i$-th refinement level is denoted as $\kappa_i$. The increase in condition numbers is as expected from the nature of the two problems – proportional to $\mathfrak{n}$ and $\mathfrak{n}^2$ for $P_1$ and $P_2$, respectively, with $\sqrt{\mathfrak{n}}$ as a stand-in for the mesh size. Note that, from Proposition 3.8 and Figure 11, the number of dofs at the $i$-th refinement level is $\mathfrak{n}_i = \mu(2^{i+1}+1)^2 + 3$.

We start from a rectangular planar geometry defined using a mesh $\mathcal{T}$ with two interior extraordinary points (valences 3 and 5); the planar control points are chosen such that the planar geometry coincides with the projection of the cylindrical roof on the $xy$-plane. To build the cylindrical roof, we first refine the planar geometry a desired number of times and then perform an $L^2$-projection to find the height of the control points. Note that since we are trying to match a given target geometry here, it is reasonable to perform this fitting for each new refinement.

Solving the Kirchhoff–Love problem on the geometry thus obtained, the solution is as shown in Figure 12(a); the deformations have been scaled up by a factor of 15 for the purpose of visualization. The vertical displacement at the mid-point of the free edges converges toward the reference Kirchhoff–Love reference solution of $w_{ref}$; the normalised error in the computed solution is shown in Figure 12(b). Note that since we did not exploit any symmetry conditions in the simulation, we display the error against the square root of a quarter of the number of degrees of freedom $\mathfrak{n}$ for the spline space defined on $\mathcal{T}$. This brings the results in line with those of [33] which assumed four-fold symmetry and plotted the error versus the number of control points per edge of their rectangular mesh. The performance of our construction is indistinguishable from the results of [33] which used standard tensor-product B-splines.

### 5.3. Surface Cahn–Hilliard problem

We now solve the fourth-order non-linear Cahn–Hilliard problem on the topologically complex surface $\Omega$ shown in Figure 13. The non-dimensional strong form of the problem is as below (see [3] for the associated weak form):

$$\frac{\partial c}{\partial t} = \nabla_\Omega \cdot (c(1-c)\nabla_\Omega(\mathbb{N}_2\mu_c - \Delta_\Omega c)) \quad \text{on } \Omega \times [0,\, T]\,,$$
$$c(\boldsymbol{x},0) = c_0(\boldsymbol{x}) \quad \text{on } \Omega\,,$$

where $\nabla_\Omega$ and $\Delta_\Omega$ are the surface gradient and Laplace–Beltrami operators, respectively, and $\mu_c := \frac{1}{3}\log\left(\frac{c}{1-c}\right) + 1 - 2c$. The unstructured mesh had 18,432 quadrilateral elements and 18,552 degrees of freedom, and we solved the equations for initial volume fraction $\bar{c} = 0.5$ and the corresponding value of $\mathbb{N}_2$ was 41.7313. The initial value of $c$, namely $c_0$, was determined by randomly perturbing $\bar{c}$, as described in [18, 37]. The results are shown in Figure 13. Steady state was reached for the configuration in 600 time-steps with the aid of an adaptive time-stepping scheme [37]. At all times, the solution coefficients were strictly between 0 and 1. Then, since the spline basis functions are non-negative and form a partition of a unity, it easily follows that the computed solution is pointwise between 0 and 1 at all times.

### 5.4. Laplace–Beltrami eigenvalue problem

For the final set of tests, we showcase the application of our splines for reconstructing a complex CAD geometry and solving an eigenvalue problem on it. The geometry, a portion of a BMW car, is freely available as a Blender® model. This model was imported into Rhinoceros® where a quadrilateral mesh consisting of 4,482

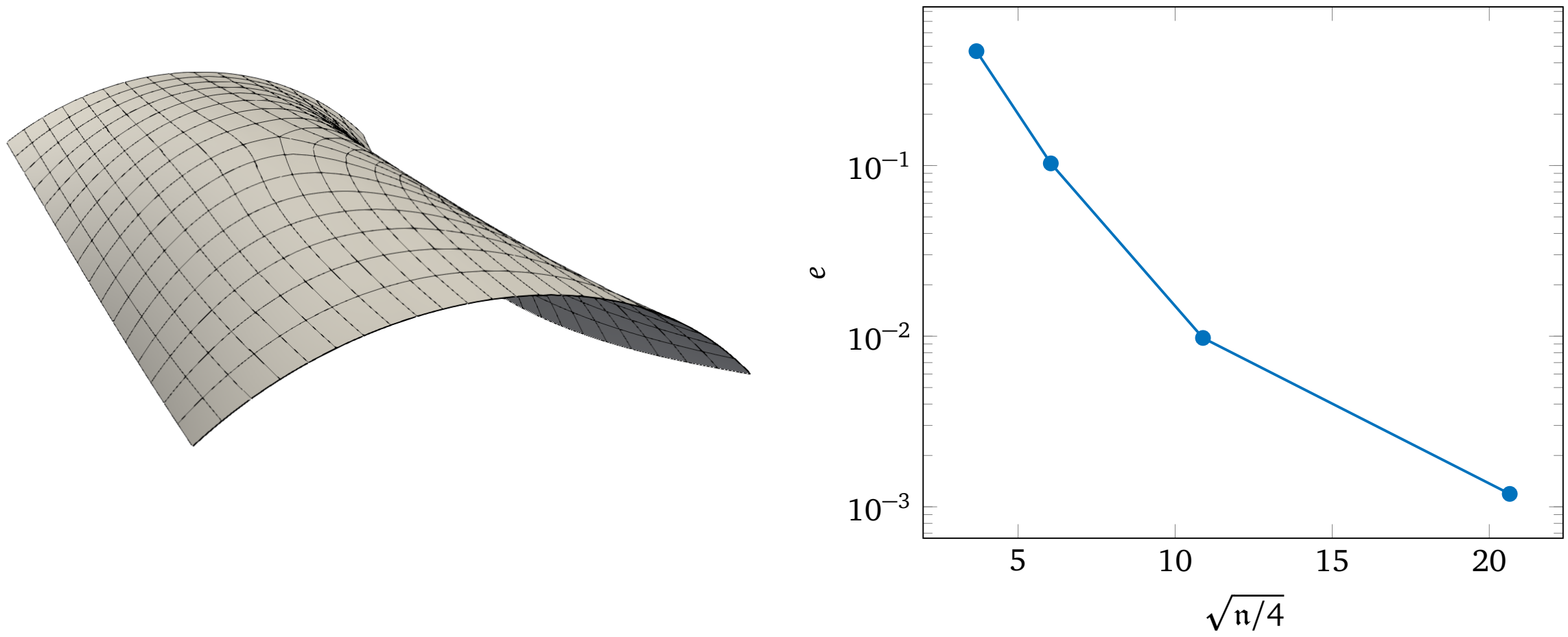

(a) deformed geometry, $\mathfrak{n} = 474$ (b) error convergence for midpoint displacement

Figure 12: For the Scordelis–Lo benchmark problem (Section 5.2), Figure (a) shows the solution at one of the refinement levels, while Figure (b) shows the convergence of the error $e := |1 - \frac{w}{w_{\text{ref}}}|$ under mesh refinement, where $w$ is the vertical displacement computed at the midpoint of one of the non-curved edges of the cylindrical roof and $w_{\text{ref}}$ is a reference solution based on Kirchhoff–Love theory [33]. Our results are based on Kirchhoff–Love thin shell theory in which transverse shear deformations are suppressed, and as one would expect, the Kirchhoff–Love converged displacement is somewhat less than that for Reissner–Mindlin theory.

faces was created; this mesh contained both boundary and interior extraordinary points. Using our splines built on this mesh (5,330 dofs), we reconstructed the geometry and solved a Laplace–Beltrami eigenvalue problem [15] on the spline geometry. This problem is defined as: find $f \in \mathscr{S}$ and $\lambda \in \mathbb{R}$ such that

$$\int_{\Omega} \nabla_{\Omega} g \cdot \nabla_{\Omega} f \; d\Omega = \int_{\Omega} g f \; d\Omega\,, \quad \forall g \in \mathscr{S}\,, \tag{30}$$

where $\mathscr{S} := \{f \in \mathscr{B} \,:\, f|_{\Gamma_b} = 0\}$ and $\Gamma_b$ is the union of the two bottom edges of the car hull. The spline geometry and select eigenmodes for the problem are shown in Figure 14.

## 6. Conclusions

We have presented the construction and refinement of almost-$C^1$ splines, that is, analysis-suitable biquadratic spline spaces on fully unstructured quadrilateral meshes, for building spline surfaces of arbitrary topology as well as for solving fourth-order problems on them. Several numerical examples of challenging fourth-order problems have been presented to exemplify this. The corresponding almost-$C^1$ spline basis functions are well-conditioned and have several B-spline-like properties such as partition of unity, non-negativity, local support and linear independence. Furthermore, we have described the construction explicitly and in a self-contained manner using Bézier extractions for enabling immediate implementation.

We use approximate smoothness in an explicit manner for our construction and strongly believe that this is a powerful approach for arriving at spline constructions that can circumvent many of the obstacles faced by strongly smooth unstructured splines, while at the same time retaining many of the latter's advantages. In this first paper we focus on the construction itself and while the numerical results are highly encouraging, a theoretical analysis of convergence is an interesting topic for future research. It is our opinion that such a theoretical analysis may be simpler for almost-$C^1$ splines than for constructions that employ singularities in the definitions of the spline basis functions (although they are strongly smooth, nested and demonstrate optimal convergence in numerical tests; the singularities usually pose significant difficulties for any standard approximation proofs). A possible way to prove convergence for almost-$C^1$ splines may be to follow strategies for non-conforming finite elements, such as in [54].

Furthermore, within this new framework, there are many extensions possible which will be the focus of our future work. Some of these are the formulation of similar constructions for higher polynomial degrees and higher orders of approximate smoothness, as well as incorporation of local refinement. For the latter, we note

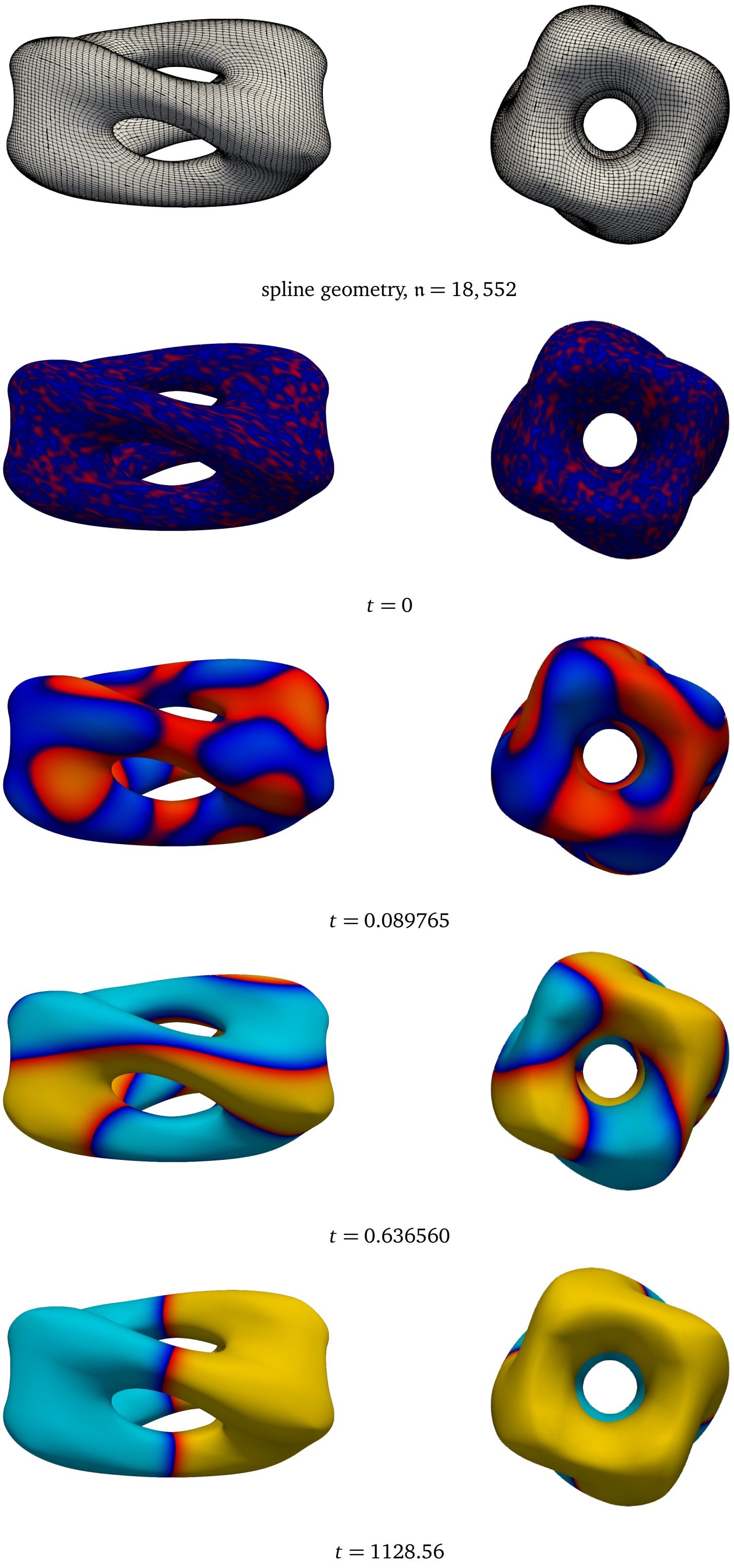

Figure 13: These figures shows the initial volume fraction distribution ($t = 0$) over a surface of non-trivial topology (top row), which is the domain of interest for the surface Cahn–Hilliard problem (Section 5.3), and the rows below show its time-evolution. The meshes used for the computation contained 18,552 degrees of freedom. (The left and right columns show different views of the geometry and solution.)

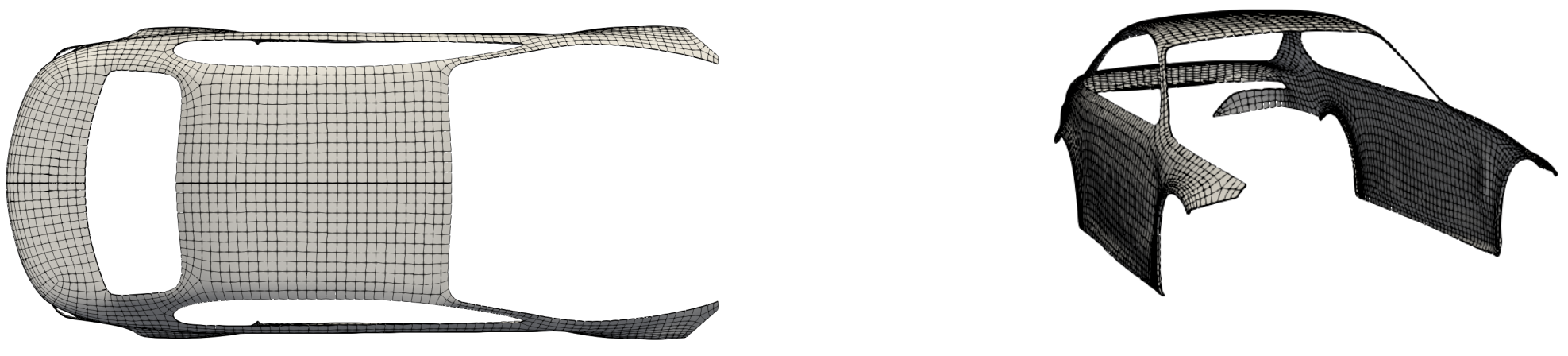

(a) spline geometry

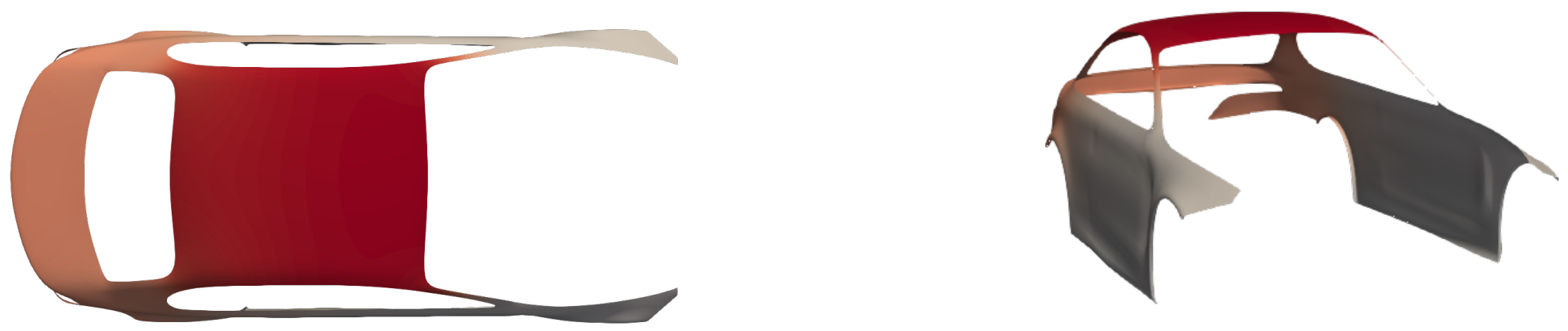

(b) mode 1

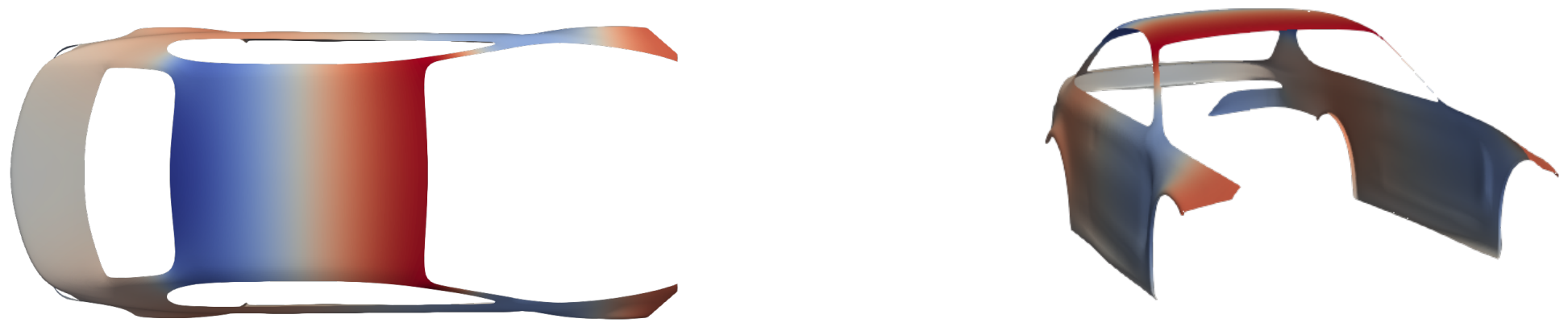

(c) mode 11

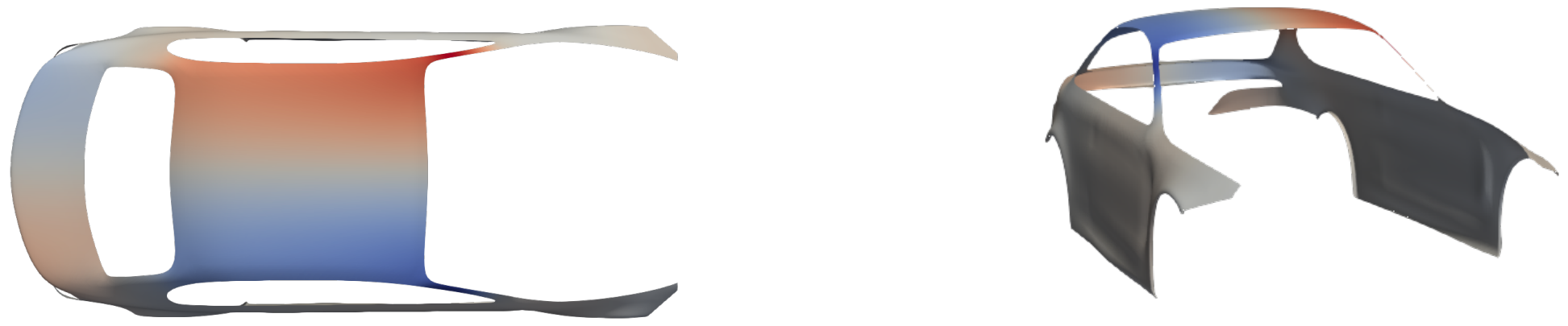

(d) mode 12

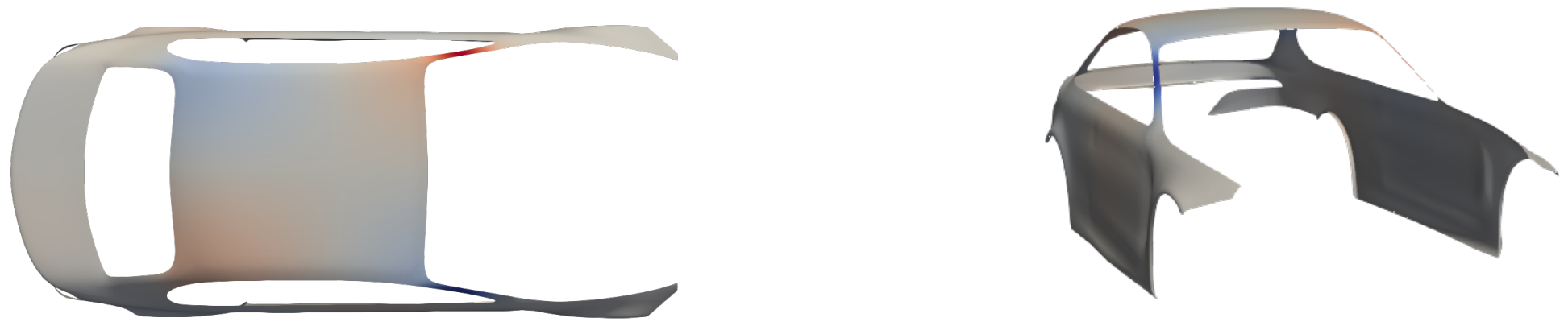

(e) mode 15

Figure 14: The figure in the top row shows a part of a BMW car – the model is a Blender® demofile, a coarse quadrilateral mesh for it was generated using Rhinoceros® and the spline geometry was constructed on that mesh using our unstructured splines. The bottom rows show select eigenmodes corresponding to the Laplace–Beltrami eigenvalue problem (Section 5.4) solved on this geometry. The mesh contains 4,482 faces and 5,330 dofs.

that since our construction is highly local, it can be readily embedded within locally refined spline constructions (e.g., similarly to [10]) by assuming sufficient separation between hanging nodes and extraordinary vertices – for biquadratic splines, a 2-ring distance between extraordinary vertices and hanging nodes is expected to be sufficient. As an alternative approach for local refinement, a generalization similar to hierarchical B-splines may also be possible, even though the spaces on different refinement levels are not nested. Finally, other minor extensions include generalizing the set of corner vertices to include boundary vertices of valence higher than 1, as well as reducing the smoothness across select interior edges to create geometries with creased features.

**Acknowledgements**

The research of Deepesh Toshniwal is supported by project number 212.150 awarded through the Veni research programme by the Dutch Research Council (NWO). The research of Thomas Takacs is partially supported by the Austrian Science Fund (FWF) and the government of Upper Austria through the project P 30926-NBL entitled "Weak and approximate $C^1$-smoothness in isogeometric analysis".

**Appendix A. Alternative construction 1: Truncation using degree elevation**

A construction similar to the one explained in Section 3.5 can be achieved if the functions $B^*_\phi \in \mathscr{B}$ are degree elevated on every extraordinary face, i.e., $B^*_\phi|_\sigma$ is represented as a bicubic polynomial for each $\sigma \in \mathscr{T}_2^E$. Then $B_\phi|_\sigma$ is given as the truncated, degree elevated version of $B^*_\phi|_\sigma$. Algorithmically, the only change is in the matrix $K$, which has to be replaced by the degree elevation matrix

$$K_{d.e.} = \begin{bmatrix} 1 & & \\ \frac{1}{3} & \frac{2}{3} & \\ & \frac{2}{3} & \frac{1}{3} \\ & & 1 \end{bmatrix} , \tag{A.1}$$

and in the local basis $b^1_{jk,\square}$ in (10), which is replaced by bicubic Bernstein polynomials. The truncation step as well as the basis construction for the extraordinary vertex splines remain the same.

**Appendix B. Alternative construction 2: Geometry-independent templates**

In the following we discuss how the basis functions $B_{(\gamma,\nu)}$, $\nu \in \{1,2,3\}$, can be defined independently of the geometry. This is achieved by replacing the projection onto a prescribed tangent plane, as in (13), by a template configuration depending only on the valence. We assume that the control triangle is always given as

$$(\boldsymbol{a}_1, \boldsymbol{a}_2, \boldsymbol{a}_3) = \left( (0,1)^T, \left(-\frac{\sqrt{3}}{2}, -\frac{1}{2}\right)^T, \left(\frac{\sqrt{3}}{2}, -\frac{1}{2}\right)^T \right) .$$

Then we define points $\boldsymbol{c}^i_{jk}$, with $(j,k) \in \{0,1\}^2$, $i \in \{1,\dots,\mu\}$ to be

$$\begin{aligned} \boldsymbol{c}^i_{00} &= (0,0)^T \\ \boldsymbol{c}^i_{11} &= \frac{1}{2}(-\sin(2\pi(i-1)/\mu), \cos(2\pi(i-1)/\mu))^T \\ \boldsymbol{c}^i_{10} = \boldsymbol{c}^{i+1}_{01} &= \frac{1}{4\cos(\pi/\mu)}(-\sin(2\pi(i-1/2)/\mu), \cos(2\pi(i-1/2)/\mu))^T , \end{aligned}$$

see Figure B.15.

Consequently, the coefficients of the basis functions $B_{(\gamma,\nu)}$ are again given as the barycentric coordinates with respect to the vertices of the triangle, i.e.,

$$B_{(\gamma,\nu)}|_{\sigma_i} = \sum_{jk=0}^{3} \hat{c}_{jk}[B_{(\gamma,\nu)};\sigma_j] b^1_{jk,\square} ,$$

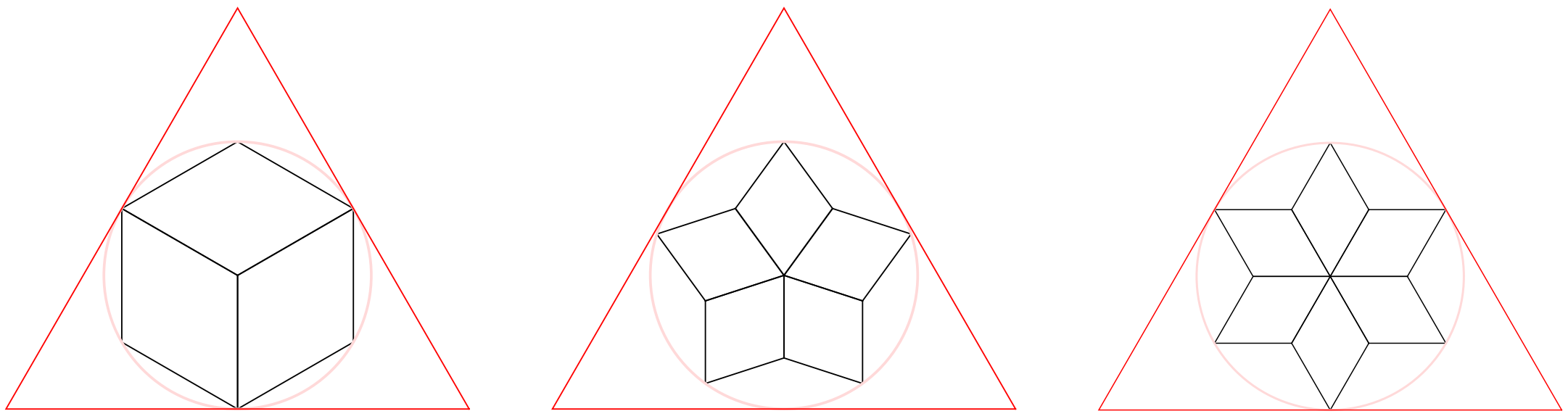

(a) valence 3. (b) valence 5. (c) valence 6.

Figure B.15: Templates for valencies 3, 5 and 6.

with

$$\hat{c}_{jk}[B_{(\gamma,\nu)};\sigma_i] = \begin{cases} \lambda_\nu(\mathbf{c}^i_{jk}), & (j,k) \in \{0,1\}^2, \\ 0, & \text{otherwise}. \end{cases}$$

Independent of the valence, we always have $\hat{c}_{00}[B_{(\gamma,\nu)};\sigma_i] = \frac{1}{3}$. For $\mu = 3$ we have

$$\hat{c}_{11}[B_{(\gamma,1)};\sigma_1] = \frac{2}{3} \qquad \hat{c}_{10}[B_{(\gamma,1)};\sigma_1] = \frac{1}{2}$$
$$\hat{c}_{11}[B_{(\gamma,1)};\sigma_2] = \frac{1}{6} \qquad \hat{c}_{10}[B_{(\gamma,1)};\sigma_2] = 0$$
$$\hat{c}_{11}[B_{(\gamma,1)};\sigma_3] = \frac{1}{6} \qquad \hat{c}_{10}[B_{(\gamma,1)};\sigma_3] = \frac{1}{2},$$

for $\mu = 5$ we have

$$\hat{c}_{11}[B_{(\gamma,1)};\sigma_1] = \frac{2}{3} \qquad \hat{c}_{11}[B_{(\gamma,2)};\sigma_1] = \frac{1}{6} \approx 0.166667$$
$$\hat{c}_{10}[B_{(\gamma,1)};\sigma_1] = \frac{1}{2} \qquad \hat{c}_{10}[B_{(\gamma,2)};\sigma_1] = \frac{3+\sqrt{15-6\sqrt{5}}}{12} \approx 0.354867$$
$$\hat{c}_{11}[B_{(\gamma,1)};\sigma_2] = \frac{3+\sqrt{5}}{12} \qquad \hat{c}_{11}[B_{(\gamma,2)};\sigma_2] = \frac{9-\sqrt{5}+\sqrt{30+6\sqrt{5}}}{24} \approx 0.556377$$
$$\hat{c}_{10}[B_{(\gamma,1)};\sigma_2] = \frac{1+\sqrt{5}}{12} \qquad \hat{c}_{10}[B_{(\gamma,2)};\sigma_2] = \frac{11-\sqrt{5}+\sqrt{30-6\sqrt{5}}}{24} \approx 0.534843$$
$$\hat{c}_{11}[B_{(\gamma,1)};\sigma_3] = \frac{3-\sqrt{5}}{12} \qquad \hat{c}_{11}[B_{(\gamma,2)};\sigma_3] = \frac{9+\sqrt{5}+\sqrt{30-6\sqrt{5}}}{24} \approx 0.637848$$
$$\hat{c}_{10}[B_{(\gamma,1)};\sigma_3] = \frac{3-\sqrt{5}}{6} \qquad \hat{c}_{10}[B_{(\gamma,2)};\sigma_3] = \frac{3+\sqrt{5}}{12} \approx 0.436339$$
$$\hat{c}_{11}[B_{(\gamma,1)};\sigma_4] = \frac{3-\sqrt{5}}{12} \qquad \hat{c}_{11}[B_{(\gamma,2)};\sigma_4] = \frac{9+\sqrt{5}-\sqrt{30-6\sqrt{5}}}{24} \approx 0.298491$$
$$\hat{c}_{10}[B_{(\gamma,1)};\sigma_4] = \frac{1+\sqrt{5}}{12} \qquad \hat{c}_{10}[B_{(\gamma,2)};\sigma_4] = \frac{11-\sqrt{5}-\sqrt{30-6\sqrt{5}}}{24} \approx 0.195485$$
$$\hat{c}_{11}[B_{(\gamma,1)};\sigma_5] = \frac{3+\sqrt{5}}{12} \qquad \hat{c}_{11}[B_{(\gamma,2)};\sigma_5] = \frac{9-\sqrt{5}-\sqrt{30+6\sqrt{5}}}{24} \approx 0.00728413$$
$$\hat{c}_{10}[B_{(\gamma,1)};\sigma_5] = \frac{1}{2} \qquad \hat{c}_{10}[B_{(\gamma,2)};\sigma_5] = \frac{3-\sqrt{15-6\sqrt{5}}}{12} \approx 0.145133,$$

and for $\mu = 6$ we have

$$\hat{c}_{11}[B_{(\gamma,1)};\sigma_1] = \frac{2}{3} \qquad \hat{c}_{10}[B_{(\gamma,1)};\sigma_1] = \frac{1}{2}$$
$$\hat{c}_{11}[B_{(\gamma,1)};\sigma_2] = \frac{1}{2} \qquad \hat{c}_{10}[B_{(\gamma,1)};\sigma_2] = \frac{1}{3}$$
$$\hat{c}_{11}[B_{(\gamma,1)};\sigma_3] = \frac{1}{6} \qquad \hat{c}_{10}[B_{(\gamma,1)};\sigma_3] = \frac{1}{6}$$
$$\hat{c}_{11}[B_{(\gamma,1)};\sigma_4] = 0 \qquad \hat{c}_{10}[B_{(\gamma,1)};\sigma_4] = \frac{1}{6}$$
$$\hat{c}_{11}[B_{(\gamma,1)};\sigma_5] = \frac{1}{6} \qquad \hat{c}_{10}[B_{(\gamma,1)};\sigma_5] = \frac{1}{3}$$
$$\hat{c}_{11}[B_{(\gamma,1)};\sigma_6] = \frac{1}{2} \qquad \hat{c}_{10}[B_{(\gamma,1)};\sigma_6] = \frac{1}{2}\,.$$

The other coefficients are determined via the symmetry of the configuration. While valencies 3 and 6 allow rotationally symmetric configurations, where all three basis functions look alike, valence 5 yields two different types of functions ($B_{(\gamma,2)}$ and $B_{(\gamma,3)}$ are equivalent up to reflection). Similar templates can be derived for extraordinary vertices at the boundary. We show three such configurations in Figure B.16 and leave the computation of coefficients as an exercise to the readers.

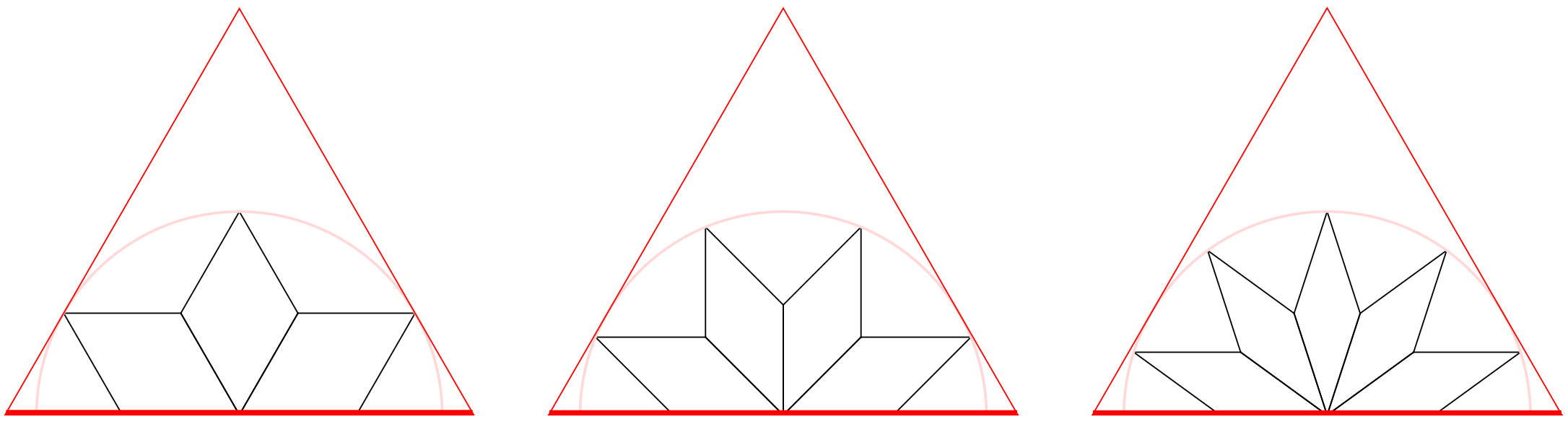

(a) valence 3. (b) valence 4. (c) valence 5.

Figure B.16: Templates for boundary vertices with valencies 3, 4 and 5. The boundary is assumed to be at the bottom edge of the control triangle.

## Appendix C. Alternative construction 3: A smooth subspace without local subdivison

In the following we define an approximate $C^1$ subspace $\mathscr{B}^\dagger$ of the B-spline space $\mathscr{B}^*$ defined in Section 3.1. For this we need to reduce the number of degrees of freedom. Instead of taking all faces in $\mathscr{T}_2$ as degrees of freedom, we mark faces in such a way, that for each extraordinary vertex exactly three faces in its 1-ring are marked. In addition, all faces that are not contained in the 1-ring of an extraordinary vertex are marked. Such a marking may not always be possible. However, by bisecting the mesh once, all 1-rings of extraordinary vertices are disjoint and a valid marking exists. We collect all marked faces in $\mathscr{T}_2^\dagger$.

Thus, we end up with the following degrees of freedom:

- Face dofs: We associate one degree of freedom to each marked face $\sigma \in \mathscr{T}_2^\dagger$.
- Boundary edge dofs: We associate one degree of freedom to each boundary edge $\tau \in \mathscr{T}_1^B$.
- Corner vertex dofs: We associate one degree of freedom to each corner vertex $\gamma \in \mathscr{T}_0^C$.

Given an extraordinary vertex $\gamma \in \mathscr{T}_0^E$ of valence $\mu$ and let $\sigma_1, \sigma_2, \ldots, \sigma_\mu$ be the faces surrounding it. Moreover, let $\{i_1, i_2, i_3\} \subset \{1, 2, \ldots, \mu\}$ denote the marked faces $\sigma_{i_1}$, $\sigma_{i_2}$ and $\sigma_{i_3}$. Given a spline geometry $\boldsymbol{x}^*$, let $\boldsymbol{x}_i^*$ be the control points corresponding to the faces $\sigma_i$.

Let again $P$ be the orthogonal projection onto the tangent plane given by a prescribed normal $\boldsymbol{n}_\gamma$. This yields new control points $P(\boldsymbol{x}_i^*)$. New basis functions $B^\dagger_{\sigma_{i_\nu}}$ can now be constructed for $\sigma_{i_1}$, $\sigma_{i_2}$ and $\sigma_{i_3}$, via

$$B^\dagger_{\sigma_{i_\nu}} = \sum_{i=1}^{\mu} \lambda_\nu(P(\boldsymbol{x}_i^*)) B^*_{\sigma_i}\,,$$

where $\lambda_\nu(\boldsymbol{c})$ denotes the $\nu$-th component of the barycentric coordinates of $\boldsymbol{c}$ with respect to the triangle

$$(P(\boldsymbol{x}^*_{i_1}), P(\boldsymbol{x}^*_{i_2}), P(\boldsymbol{x}^*_{i_3}))\ .$$

All functions corresponding to regular faces, boundary edges and corner vertices remain unchanged, i.e., $B^\dagger_\phi = B^*_\phi$ for all $\phi \in \mathscr{T}_2 \setminus \mathscr{T}_2^E \cup \mathscr{T}_1^B \cup \mathscr{T}_0^C$.

Similar to the construction in Section 3.5 one can describe a basis with respect to a control triangle which is different from the control points of marked faces. Such a construction based on a control triangle is only feasible if no face contains more than one extraordinary vertex. Alternatively, a construction based on a template as in Appendix B can be applied as well.

## References

[1] Pieter J Barendrecht. Isogeometric analysis for subdivision surfaces. *Eindhoven University of Technology: Eindhoven, The Netherlands*, 2013.

[2] Pieter J Barendrecht, Michael Bartoň, and Jiří Kosinka. Efficient quadrature rules for subdivision surfaces in isogeometric analysis. *Computer Methods in Applied Mechanics and Engineering*, 340:1–23, 2018.

[3] A. Bartezzaghi, L. Dedè, and A. Quarteroni. Isogeometric analysis of high order partial differential equations on surfaces. *Computer Methods in Applied Mechanics and Engineering*, 295:446–469, 2015.

[4] Andrea Benvenuti and Giancarlo Sangalli. *Isogeometric Analysis for $C^1$-continuous Mortar Method*. PhD thesis, PhD thesis, Corso di Dottorato in Matematica e Statistica, Università degli studi di Pavia, 2017.

[5] Michel Bercovier and Tanya Matskewich. *Smooth Bézier Surfaces over Unstructured Quadrilateral Meshes*. Lecture Notes of the Unione Matematica Italiana. Springer International Publishing, 2017.

[6] Paul T Boggs, Alan Althsuler, Alex R Larzelere, Edward J Walsh, Ruuobert L Clay, and Michael F Hardwick. Dart system analysis. Technical report, Sandia National Laboratories, 2005.

[7] M. J. Borden, M. A. Scott, J. A. Evans, and T. J. R. Hughes. Isogeometric finite element data structures based on Bézier extraction of NURBS. *International Journal for Numerical Methods in Engineering*, 87:15–47, 2011.

[8] D. Burkhart, B. Hamann, and G. Umlauf. Iso-geometric Finite Element Analysis Based on Catmull-Clark Subdivision Solids. 29(5):1575–1584, 2010.

[9] Hugo Casquero, Carles Bona-Casas, Deepesh Toshniwal, Thomas J.R. Hughes, Hector Gomez, and Yongjie Jessica Zhang. The divergence-conforming immersed boundary method: Application to vesicle and capsule dynamics. *Journal of Computational Physics*, 425:109872, 2021.

[10] Hugo Casquero, Xiaodong Wei, Deepesh Toshniwal, Angran Li, T.J.R. Hughes, Josef Kiendl, and Yongjie Jessica Zhang. Seamless integration of design and Kirchhoff–Love shell analysis using analysis-suitable unstructured T-splines. *Computer Methods in Applied Mechanics and Engineering*, 360:112765, 2020.

[11] Edwin Catmull and James Clark. Recursively generated b-spline surfaces on arbitrary topological meshes. *Computer-aided design*, 10(6):350–355, 1978.

[12] Chiu Ling Chan, Cosmin Anitescu, and Timon Rabczuk. Isogeometric analysis with strong multipatch $C^1$-coupling. *Computer Aided Geometric Design*, 62:294–310, 2018.

[13] Chiu Ling Chan, Cosmin Anitescu, and Timon Rabczuk. Strong multipatch $C^1$-coupling for isogeometric analysis on 2D and 3D domains. *Computer Methods in Applied Mechanics and Engineering*, 357:112599, 2019.

[14] A. Collin, G. Sangalli, and T. Takacs. Analysis-suitable $G^1$ multi-patch parametrizations for $C^1$ isogeometric spaces. *Computer Aided Geometric Design*, 47:93–113, 2016.

[15] Luca Dedè and Alfio Quarteroni. Isogeometric analysis for second order partial differential equations on surfaces. *Computer Methods in Applied Mechanics and Engineering*, 284:807–834, 2015.

[16] Daniel Doo and Malcolm Sabin. Behaviour of recursive division surfaces near extraordinary points. *Computer-Aided Design*, 10(6):356–360, 1978.

[17] G. E. Farin, J. Hoschek, and M.-S. Kim. *Handbook of Computer Aided Geometric Design*. Elsevier, 2002.

[18] H. Gómez, V. M. Calo, Y. Bazilevs, and T. J. R. Hughes. Isogeometric analysis of the Cahn–Hilliard phase-field model. *Computer Methods in Applied Mechanics and Engineering*, 197:4333–4352, 2008.

[19] C. M. Grimm and J. F. Hughes. Modeling surfaces of arbitrary topology using manifolds. In *Proceedings of the 22nd Annual Conference on Computer Graphics and Interactive Techniques*, pages 359–368. ACM Press, 1995.

[20] Jan Grošelj, Mario Kapl, Marjeta Knez, Thomas Takacs, and Vito Vitrih. A super-smooth c1 spline space over planar mixed triangle and quadrilateral meshes. *Computers & Mathematics with Applications*, 80(12):2623 – 2643, 2020.

[21] Yujie Guo and Martin Ruess. Nitsche's method for a coupling of isogeometric thin shells and blended shell structures. *Computer Methods in Applied Mechanics and Engineering*, 284:881–905, 2015.

[22] Thomas Horger, Alessandro Reali, Barbara Wohlmuth, and Linus Wunderlich. A hybrid isogeometric approach on multi-patches with applications to Kirchhoff plates and eigenvalue problems. *Computer Methods in Applied Mechanics and Engineering*, 348:396–408, 2019.

[23] T. J. R. Hughes. *The Finite Element Method: Linear Static and Dynamic Finite Element Analysis*. Courier Corporation, 2012.

[24] T. J. R. Hughes, J. A. Cottrell, and Y. Bazilevs. Isogeometric analysis: CAD, finite elements, NURBS, exact geometry and mesh refinement. *Computer Methods in Applied Mechanics and Engineering*, 194:4135–4195, 2005.

[25] Thomas J.R. Hughes, Giancarlo Sangalli, Thomas Takacs, and Deepesh Toshniwal. Smooth multi-patch discretizations in isogeometric analysis. Handbook of Numerical Analysis. Elsevier, 2020.

[26] M. Kapl and V. Vitrih. Space of $C^2$-smooth geometrically continuous isogeometric functions on two-patch geometries. *Computers & Mathematics with Applications*, 73:37–59, 2017.

[27] M. Kapl, V. Vitrih, B. Jüttler, and K. Birner. Isogeometric analysis with geometrically continuous functions on two-patch geometries. *Computers & Mathematics with Applications*, 70:1518–1538, 2015.

[28] Mario Kapl, Giancarlo Sangalli, and Thomas Takacs. Dimension and basis construction for analysis-suitable $G^1$ two-patch parameterizations. *Computer Aided Geometric Design*, 52-53:75–89, 2017.
[29] Mario Kapl, Giancarlo Sangalli, and Thomas Takacs. Construction of analysis-suitable $G^1$ planar multi-patch parameterizations. *Computer-Aided Design*, 97:41–55, 2018.
[30] Mario Kapl, Giancarlo Sangalli, and Thomas Takacs. An isogeometric $C^1$ subspace on unstructured multi-patch planar domains. *Computer Aided Geometric Design*, 69:55–75, 2019.
[31] Mario Kapl, Giancarlo Sangalli, and Thomas Takacs. A family of $C^1$ quadrilateral finite elements. *Advances in Computational Mathematics*, 47(6):1–38, 2021.
[32] Mario Kapl and Vito Vitrih. Dimension and basis construction for $C^2$-smooth isogeometric spline spaces over bilinear-like $G^2$ two-patch parameterizations. *Journal of Computational and Applied Mathematics*, 335:289–311, 2018.
[33] Josef Kiendl, K-U Bletzinger, Johannes Linhard, and Roland Wüchner. Isogeometric shell analysis with kirchhoff–love elements. *Computer methods in applied mechanics and engineering*, 198(49-52):3902–3914, 2009.
[34] Kim Jie Koh, Deepesh Toshniwal, and Fehmi Cirak. An optimally convergent smooth blended b-spline construction for unstructured quadrilateral and hexahedral meshes. *arXiv preprint arXiv:2111.04401*, 2021.
[35] Ming-Jun Lai and Larry L Schumaker. *Spline functions on triangulations*, volume 110. Cambridge University Press, 2007.
[36] Xin Li, Xiaodong Wei, and Yongjie Jessica Zhang. Hybrid non-uniform recursive subdivision with improved convergence rates. *Computer Methods in Applied Mechanics and Engineering*, 352:606–624, 2019.
[37] J. Liu, L. Dedè, J. A. Evans, M. J. Borden, and T. J. R. Hughes. Isogeometric analysis of the advective Cahn–Hilliard equation: Spinodal decomposition under shear flow. *Journal of Computational Physics*, 242:321–350, 2013.
[38] M. Majeed and F. Cirak. Isogeometric analysis using manifold-based smooth basis functions. *Computer Methods in Applied Mechanics and Engineering*, 316:547–567, 2017.
[39] Di Miao, Zhihui Zou, Michael A Scott, Michael J Borden, and Derek C Thomas. Isogeometric Bézier dual mortaring: The Kirchhoff–Love shell problem. *Computer Methods in Applied Mechanics and Engineering*, 382:113873, 2021.
[40] Stephen Edward Moore. Discontinuous Galerkin Isogeometric Analysis for the biharmonic equation. *Computers & Mathematics with Applications*, 76(4):673–685, 2018.
[41] Bernard Mourrain, Raimundas Vidunas, and Nelly Villamizar. Dimension and bases for geometrically continuous splines on surfaces of arbitrary topology. *Computer Aided Geometric Design*, 45:108–133, 2016.
[42] T. Nguyen, K. Karčiauskas, and J. Peters. A comparative study of several classical, discrete differential and isogeometric methods for solving Poisson's equation on the disk. *Axioms*, 3:280–299, 2014.
[43] T. Nguyen and J. Peters. Refinable $C^1$ spline elements for irregular quad layout. *Computer Aided Geometric Design*, 43:123–130, 2016.
[44] United States Department of Transportation. National Highway Traffic Safety Administration. `https://www.nhtsa.gov/crash-simulation-vehicle-models`. [Online; accessed 01-December-2021].
[45] J. Peters and U. Reif. *Subdivision Surfaces*. Springer-Verlag, 2008.
[46] Katharina Rafetseder and Walter Zulehner. A new mixed approach to Kirchhoff–Love shells. *Computer Methods in Applied Mechanics and Engineering*, 346:440–455, 2019.
[47] U. Reif. A refineable space of smooth spline surfaces of arbitrary topological genus. *Journal of Approximation Theory*, 90:174–199, 1997.
[48] U. Reif. TURBS-topologically unrestricted rational B-splines. *Constructive Approximation*, 14:57–77, 1998.
[49] Ulrich Reif. Biquadratic G-spline surfaces. *Computer Aided Geometric Design*, 12(2):193–205, 1995.
[50] Andreas Riffnaller-Schiefer, Ursula H Augsdörfer, and Dieter W Fellner. Isogeometric shell analysis with NURBS compatible subdivision surfaces. *Applied Mathematics and Computation*, 272:139–147, 2016.
[51] M. A. Scott, M. J. Borden, C. V. Verhoosel, T. W. Sederberg, and T. J. R. Hughes. Isogeometric finite element data structures based on Bézier extraction of T-splines. *International Journal for Numerical Methods in Engineering*, 88:126–156, 2011.
[52] M. A. Scott, R. N. Simpson, J. A. Evans, S. Lipton, S. P. A. Bordas, T. J. R. Hughes, and T. W. Sederberg. Isogeometric boundary element analysis using unstructured T-splines. *Computer Methods in Applied Mechanics and Engineering*, 254:197–221, 2013.
[53] Kendrick M Shepherd, Xianfeng David Gu, and Thomas JR Hughes. Isogeometric model reconstruction of open shells via ricci flow and quadrilateral layout-inducing energies. *Engineering Structures*, 252:113602, 2022.
[54] Zhong Ci Shi. The fem test for convergence of nonconforming finite elements. *Mathematics of computation*, 49(180):391–405, 1987.
[55] H. Speleers. Algorithm 999: Computation of multi-degree B-splines. *ACM Transactions on Mathematical Software*, 45:Article No. 43, 2019.
[56] Hendrik Speleers. A normalized basis for reduced clough–tocher splines. *Computer aided geometric design*, 27(9):700–712, 2010.
[57] Hendrik Speleers and Deepesh Toshniwal. A general class of c1 smooth rational splines: Application to construction of exact ellipses and ellipsoids. *Computer-Aided Design*, 132:102982, 2021.
[58] J. Stam. Exact evaluation of Catmull–Clark subdivision surfaces at arbitrary parameter values. In *Proceedings of the 25th Annual Conference on Computer Graphics and Interactive Techniques*, pages 395–404. ACM Press, 1998.
[59] Thomas Takacs. Construction of smooth isogeometric function spaces on singularly parameterized domains. In *International Conference on Curves and Surfaces*, pages 433–451. Springer, 2014.
[60] Thomas Takacs and Bert Jüttler. $H^2$ regularity properties of singular parameterizations in isogeometric analysis. *Graphical models*, 74(6):361–372, 2012.
[61] D. Toshniwal, H. Speleers, R. R. Hiemstra, and T. J. R. Hughes. Multi-degree smooth polar splines: A framework for geometric modeling and isogeometric analysis. *Computer Methods in Applied Mechanics and Engineering*, 316:1005–1061, 2017.
[62] D. Toshniwal, H. Speleers, R. R. Hiemstra, C. Manni, and T. J. R. Hughes. Multi-degree B-splines: Algorithmic computation and properties. *Computer Aided Geometric Design*, 76:Article No. 101792, 2020.
[63] D. Toshniwal, H. Speleers, and T. J. R. Hughes. Smooth cubic spline spaces on unstructured quadrilateral meshes with particular emphasis on extraordinary points: Geometric design and isogeometric analysis considerations. *Computer Methods in Applied Mechanics and Engineering*, 327:411–458, 2017.
[64] Deepesh Toshniwal. Quadratic splines on quad-tri meshes: Construction and an application to simulations on watertight reconstructions of trimmed surfaces. *Computer Methods in Applied Mechanics and Engineering*, 388:114174, 2022.
[65] Deepesh Toshniwal and Thomas J.R. Hughes. Isogeometric discrete differential forms: Non-uniform degrees, bézier extraction, polar splines and flows on surfaces. *Computer Methods in Applied Mechanics and Engineering*, 376:113576, 2021.
[66] X. Wei, Y. Zhang, T. J. R. Hughes, and M. A. Scott. Truncated hierarchical Catmull–Clark subdivision with local refinement. *Computer*

*Methods in Applied Mechanics and Engineering*, 291:1–20, 2015.

[67] X. Wei, Y. J. Zhang, D. Toshniwal, H. Speleers, X. Li, C. Manni, J. A Evans, and T. J. R. Hughes. Blended B-spline construction on unstructured quadrilateral and hexahedral meshes with optimal convergence rates in isogeometric analysis. *Computer Methods in Applied Mechanics and Engineering*, 341:609–639, 2018.

[68] Pascal Weinmüller and Thomas Takacs. Construction of approximate $C^1$ bases for isogeometric analysis on two-patch domains. *Computer Methods in Applied Mechanics and Engineering*, 385:114017, 2021.

[69] Qiaoling Zhang, Malcolm Sabin, and Fehmi Cirak. Subdivision surfaces with isogeometric analysis adapted refinement weights. *Computer-Aided Design*, 102:104–114, 2018.

[70] Qiaoling Zhang, Thomas Takacs, and Fehmi Cirak. Manifold-based b-splines on unstructured meshes. In *Conference on Isogeometric Analysis and Applications*, pages 243–262. Springer, 2018.

[71] Christopher Zimmermann, Deepesh Toshniwal, Chad M. Landis, Thomas J. R. Hughes, Kranthi K. Mandadapu, and Roger A. Sauer. An isogeometric finite element formulation for phase transitions on deforming surfaces. 2019.